\documentclass[a4paper]{amsart}
\usepackage{amssymb,verbatim,xypic}

\def\xcirc{\objectmargin{0.1pc}\def\objectstyle{\sssize}\diagram\squarify<1pt>{}\circled\enddiagram}

\newtheorem{theorem}{Theorem}[section]
\newtheorem{lemma}[theorem]{Lemma}
\newtheorem{proposition}[theorem]{Proposition}

\theoremstyle{definition}

\theoremstyle{remark}
\newtheorem{remark}[theorem]{Remark}

\DeclareMathOperator{\Der}{Der}

\DeclareMathOperator{\HS}{H}

\DeclareMathOperator{\HH}{HH}

\DeclareMathOperator{\ide}{id}

\DeclareMathOperator{\Hom}{Hom}

\DeclareMathOperator{\bN}{\mathbb{N}}

\DeclareMathOperator{\bQ}{\mathbb{Q}}

\DeclareMathOperator{\sg}{sg}

\newcommand{\ot}{\otimes}
\newcommand{\wt}{\widetilde}
\newcommand{\wh}{\widehat}
\newcommand{\ov}{\overline}
\newcommand{\ep}{\epsilon}
\newcommand{\de}{\delta}
\newcommand{\la}{\lambda}
\newcommand{\si}{\sigma}
\newcommand{\im}{\imath}
\newcommand{\jm}{\jmath}
\newcommand{\fg}{\mathfrak{g}}
\newcommand{\ba}{\mathbf a}
\newcommand{\bc}{\mathbf c}
\newcommand{\bx}{\mathbf x}
\newcommand{\bv}{\mathbf v}

\begin{document}

\title{Cohomology ring of differential operator rings}


\author{Graciela Carboni}
\address{Departamento de Matem\'atica\\ Facultad de Ciencias Exactas y Naturales, Pabell\'on
1 - Ciudad Universitaria\\ (1428) Buenos Aires, Argentina.} \curraddr{}\email{}
\thanks{Supported UBACYT 095}

\author{Jorge A. Guccione}
\address{Departamento de Matem\'atica\\ Facultad de Ciencias Exactas y Naturales, Pabell\'on
1 - Ciudad Universitaria\\ (1428) Buenos Aires, Argentina.} \curraddr{}\email{vander@dm.uba.ar}
\thanks{UBACYT 095 and PIP 112-200801-00900 (CONICET)}

\author{Juan J. Guccione}
\address{Departamento de Matem\'atica\\ Facultad de Ciencias Exactas y Naturales\\
Pabell\'on 1 - Ciudad Universitaria\\ (1428) Buenos Aires,
Argentina.}\curraddr{}\email{jjgucci@dm.uba.ar}
\thanks{UBACYT 095 and PIP 112-200801-00900 (CONICET)}


\date{}

\dedicatory{}


\begin{abstract} We compute the multiplicative structure in the Hocshchild cohomology ring of a differential operators ring and the cap product of Hochschild cohomology on the Hochschild homology.
\end{abstract}

\maketitle

\section*{Introduction}
Let $k$ be a field and $A$ an associative $k$-algebra with $1$. An extension $E/A$ of $A$ is a {\em differential operator ring} on $A$ if there exists a Lie $k$-algebra $\fg$ and a $k$-vector space embedding $x\mapsto\ov{x}$, of $\fg$ into $E$, such that for all $x,y\in\fg$, $a\in A$:

\smallskip

\begin{enumerate}

\item $\ov{x}a-a\ov{x}= a^x$, where $a\mapsto a^x$ is a derivation,

\smallskip

\item $\ov{x}\ov{y}-\ov{y}\ov{x} =\overline{[x,y]_{\fg}} + f(x,y)$, where $[-,-]_{\fg}$ is the bracket of $\fg$ and $f\colon\fg\times \fg\to A$ is a $k$-bilinear map,

\smallskip

\item for a given basis $(x_i)_{i\in I}$ of $\fg$, the algebra $E$ is a free left $A$-module with the standard monomials in the $x_i$'s as a basis.

\end{enumerate}

\smallskip

This general construction was introduced in~\cite{Ch} and~\cite{Mc-R}. Several particular cases of this type of extensions have been considered previously in the literature. For instance:

\begin{itemize}

\smallskip

\item[-] when $\fg$ is one dimensional and $f$ is trivial, $E$ is the Ore extension $A[x,\delta]$, where $\delta(a)= a^x$,

\smallskip

\item[-] when $A=k$, one obtain the algebras studied by Sridharan in~\cite{S}, which are the quasi-commutative algebras $E$, whose associated graded algebra is a symmetric algebra,

\smallskip

\item[-] in~\cite[\S2]{Mc} this type of extensions was studied under the hypothesis that $A$ is commutative and $(x,a)\mapsto a^x$ is an action, and in~\cite[Theorem~4.2]{B-G-R} the case in which the cocycle is trivial was considered.

\end{itemize}

\medskip

In~\cite{B-C-M} and~\cite{D-T} the study of the crossed products $A\#_f H$ of a $k$-algebra $A$ by a Hopf $k$-algebra $H$ was begun, and in~\cite{M} was proved that the differential operator rings on $A$ are the crossed products of $A$ by enveloping algebras of Lie algebras.

\smallskip

In~\cite{G-G1} complexes, simpler than the canonical ones, giving the Hochschild homology and cohomology of a differential operator ring $E$ with coefficients in an $E$-bimodule $M$, were obtained. In this paper we continue this investigation by studying the Hocshchild cohomology ring of $E$ and the cap product
$$
\HS_p(E,M)\times \HH^q(E) \to \HS_{p-q}(E,M)\qquad (q\le p),
$$
in terms of the above mentioned complexes. Moreover we generalize the results of~\cite{G-G1} by considering the (co)homology of $E$ relative to a subalgebra $K$ of $A$ which is stable under the action of $\fg$ (which we also call the Hochschild (co)homology of the $K$-algebra $E$), and we seized the opportunity to fix some minor mistakes and to simplify some proofs in~\cite{G-G1}.

\smallskip

The paper is organized in the following way: In Section~1 we obtain a projective resolution $(X_*,d_*)$ of the $E$-bimodule $E$, relative to the family of all epimorphism of $E$-bimodules which split as $(E,K)$-bimodule maps. In Section~2 we determine and study comparison maps
between $(X_*,d_*)$ and the relative to $K$ normalized Hochschild resolution $(E\ot_K\ov{E}^{ \ot_K^*}\ot_K E,b'_*)$ of $E$. In Sections~3 and~4 we apply the above results in order to obtain complexes $(\ov{X}^K_*(M),\ov{d}_*)$ and $(\ov{X}_K^*(M),\ov{d}^*)$, simpler that the canonical ones, giving Hochschid homology and cohomology of the $K$-algebra $E$ with coefficients in an $E$-bimodule $M$, respectively. The main results are Theorems~\ref{th 3.4} and~\ref{th 4.4}, in which we obtain morphisms
$$
\ov{X}_K^*(E)\otimes \ov{X}_K^*(E)\to \ov{X}_K^*(E)\quad\text{and}\quad \ov{X}^K_*(M)\otimes \ov{X}_K^*(E)\to \ov{X}^K_*(M),
$$
inducing the cup and cap product, respectively. Finally in Section~5, assuming that $A$ is a symmetric algebra, we obtain further simplifications.

\section{Preliminaries}
Let $k$ be a field. In this paper all the algebras are over $k$. Let $A$ be an algebra and $H$ a Hopf algebra. We are going use the Sweedler notation $\Delta(h) = \sum_{(h)} h^{(1)}\ot_k h^{(2)}$ for the comultiplication $\Delta$ of $H$. A {\em weak action} of $H$ on $A$ is a $k$-bilinear map $(h,a)\mapsto a^h$, from $H\times A$ to $A$, such that

\smallskip

\begin{enumerate}

\item $(ab)^h =\sum_{(h)} a^{h^{(1)}}b^{h^{(2)}}$,

\smallskip

\item $1^h =\ep(h)1$,

\smallskip

\item $a^1 = a$,

\end{enumerate}

\smallskip

\noindent for $h\in H$, $a,b\in A$. By an {\em action} of $H$ on $A$ we mean a weak action such that
$$
(a^l)^h = a^{hl}\quad\text{for all $h,l\in H$, $a\in A$}.
$$

\medskip

Let $A$ be an algebra and let $H$ be a Hopf algebra acting weakly on $A$. Given a $k$-linear map $f\colon H\ot_k H\to A$ we let $A\#_f H$ denote the algebra (in general non associative and without $1$) whose underlying vector space is $A\ot_k H$ and whose multiplication is given by
$$
(a\ot_k h)(b\ot_k l) =\sum_{(h)(l)} a b^{h^{(1)}}f(h^{(2)},l^{(1)})\ot_k h^{(3)}l^{(2)},
$$
for all $a,b\in A$, $h,l\in H$. The element $a\ot_k h$ of $A\#_f H$ will usually be written $a\# h$. The algebra $A\#_f H$ is called a {\em crossed product} if it is associative with $1\# 1$ as identity element. In~\cite{B-C-M} it was proven that this happens if and only if the map $f$ and the weak action of $H$ on $A$ satisfy the following conditions

\smallskip

\begin{enumerate}

\item (Normality of $f$) for all $h\in H$ we have $f(h,1)=f(1,h) =\ep(h)1_A$,

\smallskip

\item (Cocycle condition) for all $h,l,m\in H$ we have
$$
\qquad\sum_{(h)(l)(m)}f\bigl(l^{(1)},m^{(1)}\bigr)^{h^{(1)}} f\bigl( h^{(2)},l^{(2)}m^{(2)} \bigr) = \sum_{(h)(l)} f\bigl(h^{(1)}, l^{(1)}\bigr) f\bigl(h^{(2)}l^{(2)},m\bigr),
$$

\smallskip

\item (Twisted module condition) for all $h,l\in H$ and $a\in A$ we have
$$
\sum_{(h)(l)}\bigl(a^{l^{(1)}}\bigr)^{h^{(1)}} f\bigl(h^{(2)},l^{(2)}\bigr) =\sum_{(h)(l)} f\bigl(h^{(1)},l^{(1)}\bigr)a^{h^{(2)}l^{(2)}}.
$$

\end{enumerate}

\smallskip

From now on we assume that $H$ is the enveloping algebra $U(\fg)$ of a Lie algebra~$\fg$. In this case, item~(1) of the definition of weak action implies that
$$
(ab)^x = a^xb+ab^x
$$
for each $x\in\fg$ and $a,b\in A$. So, a weak action determines a $k$-linear map
$$
\delta\colon\fg\rightarrow\Der_k(A)
$$
by $\delta(x)(a) = a^x$. Moreover if $(h,a)\mapsto a^h$ is an action, then $\de$ is a homomorphism of Lie algebras. Conversely, given a $k$-linear map $\de\colon\fg \rightarrow \Der_k(A)$, there exists a (generality non-unique) weak action of $U(\fg)$ on $A$ such that $\de(x)(a) = a^x$. When $\delta$ is a homomorphism of Lie algebras, there is a unique action of $U(\fg)$ on $A$ such that $\de(x)(a) = a^x$. For a proof of these facts see~\cite{B-C-M}. It is easy to see that each normal cocycle
$$
f\colon U(\fg)\ot_k U(\fg)\to A
$$
is convolution invertible. For a proof see~\cite[Remark 1.1]{G-G1}.

\smallskip

Next we recall some results and notations from~\cite{G-G1} that we will need later. Let $K$ be a subalgebra of $A$ which is stable under the weak action of $\fg$ (that is $\lambda^x\in K$ for all $\lambda\in K$ and $x\in \fg$) and let $E = A\#_f U(\fg)$ be a crossed product. We are going to modify the sign of some boundary maps in order to obtain simple expressions for the comparison maps.

\smallskip

To begin, we fix some notations:

\smallskip

\begin{enumerate}

\item The unadorned tensor product $\ot$ means the tensor product $\ot_K$ over $K$,

\smallskip

\item For $B=A$ or $B=E$ and each $r\in \bN$, we write $\ov{B} = B/K$,
$$
B^r = B\ot\cdots\ot B\text{ ($r$ times) and }\, \ov{B}^r =\ov{B}\ot\cdots\ot\ov{B}\text{ ($r$ times).}
$$
Moreover, for $b\in B$ we also let $b$ denote the class of $b$ in $\ov{B}$.

\smallskip

\item For each Lie algebra $\fg$ and $s\in\bN$, we write $\fg^{\land s} = \fg\land\cdots \land \fg$ ($s$ times).

\smallskip

\item Throughout this paper we will write $\ba_{1r}$ for $a_1\ot\cdots\ot a_r\in A^r$ and $\bx_{1s}$ for  $x_1\land\cdots\land x_s\in\fg^{\land s}$

\smallskip

\item For $\ba_{1r}$ and $0\le i<j\le r$, we write $\ba_{ij} = a_i\ot\cdots\ot a_j$.

\smallskip

\item For $\bx_{1s}$ and $1\le i\le s$, we write $\bx_{1\wh{\im}s} = x_1\land\cdots \land \wh{x_i}\land\cdots\land x_s$.

\smallskip

\item For $\bx_{1s}$ and $1\le i<j\le s$,  we write $\bx_{1\wh{\im}\wh{\jm}s} = x_1\land\cdots\land\wh{x_i}\land\cdots\land\wh{x_j}\land\cdots\land x_s$.

\end{enumerate}

\smallskip

Let $Y_*$ be the graded algebra generated by $A$ and the elements $y_x$, $z_x$ ($x\in\fg$) in degree zero, the elements $e_x$ ($x\in\fg$) in degree one, and the relations
\begin{align*}
& \begin{aligned}
&y_{\la x+x'} =\la y_x + y_{x'},\\ &z_{\la x+x'} =\la z_x + z_{x'},\\ &e_{\la x+x'} =\la e_x + e_{x'},
\end{aligned}
\qquad
\begin{aligned}
&y_x a = a^x + ay_x,\\ &z_x a = a^x + az_x,\\ &e_x a = a e_x,
\end{aligned}
\qquad
\begin{aligned}
&e_{x'} y_x = y_x e_{x'} + e_{[x',x]_{\fg}},\\ &e_{x'} z_x = z_x e_{x'},\\ &e_x^2 = 0,
\end{aligned}\\
&y_{x'}y_x = y_xy_{x'} + y_{[x',x]_{\fg}} + f(x',x) - f(x,x'),\\
&z_{x'}y_x = y_x z_{x'} + z_{[x',x]_{\fg}} + f(x',x) - f(x,x'),\\
&z_{x'}z_x = z_x z_{x'} + z_{[x',x]_{\fg}} + f(x',x) - f(x,x'),
\end{align*}
where $[x',x]_{\fg}$ denotes the Lie bracket of $x'$ and $x$ in $\fg$. Note that $E$ is a subalgebra of $Y_*$ via the embedding that takes $a\in A$ to $a$ and $1\# x$ to $y_x$ for all $x\in\fg$. This gives rise to an structure of left $E$-module on $Y_*$. Similarly we consider $Y_*$ as a right $E$-module via the embedding of $E$ in $Y_*$ that takes $a\in A$ to $a$ and $1\# x$ to $z_x$ for all $x\in\fg$.

\smallskip

Let $(g_i)_{i\in I}$ be a basis of $\fg$ with indexes running on an ordered set $I$.  For each $i\in I$ let us write $y_i = y_{g_i}$, $z_i = z_{g_i}$, $e_i = e_{g_i}$ and $\rho_i = z_i-y_i$.

\begin{theorem}\label{th 1.1} Each $Y_s$ is a free left $E$-module with basis
$$
\rho_{i_1}^{m_1}e_{i_1}^{\de_1}\cdots \rho_{i_l}^{m_l} e_{i_l}^{\de_l} \qquad\left(l\ge 0\text{, }i_1<\dots<i_l\in I\text{, } m_j\ge 0\text{, }\de_j \in\{0,1\}\atop m_j+\de_j>0\text{, } \de_1+\cdots+\de_l = s\right).
$$
\end{theorem}

\begin{proof} Let $\wt{\fg}$ be the direct sum of two copies $\{y_x:x\in\fg\}$ and $\{z_x:x\in\fg\}$ of $\fg$, endowed with the bracket given by
$$
[y_{x'},y_x]_{\wt{\fg}} = y_{[x',x]_{\fg}}\quad\text{and}\quad [z_{x'},z_x]_{\wt{\fg}} = [z_{x'},y_x]_{\wt{\fg}} =  z_{[x',x]_{\fg}}.
$$
Note that $\wt{\fg}$ is the semi-direct sum arising from the adjoint action of $\fg$ on itself. Let $\pi\colon U(\wt{\fg}) \to U(\fg)$ be the algebra map defined by $\pi(y_x) = \pi(z_x) = x$. Let $\Lambda(\fg)$ be the exterior algebra generated by $\fg$. That is, the algebra generated by the elements $e_x$, with $x\in \fg$, and the relations $e_{\la x+ x'} = \la e_x + e_{x'}$ and $e_x^2 = 0$, with $\la\in k$ and $x,x'\in\fg$. Let us consider the action of $U(\wt{\fg})$ on $\Lambda(\fg)$ determined by $e_{x'}^{y_x} = e_{[x,x']_{\fg}}$ and $e_{x'}^{z_x} = 0$. The enveloping algebra $U(\wt{\fg})$ of $\wt{\fg}$ acts weakly on $A\ot_k \Lambda(\fg)$ via
$$
(a\ot_k e)^u = \sum_{(u)} a^{\pi(u^{(1)})}\ot_k e^{u^{(2)}}\quad\text{($a\in A$, $e\in \Lambda(\fg)$ and $u\in U(\wt{\fg}))$.}
$$
Moreover, the map
$$
\xymatrix{\wt{f}\colon U(\wt{\fg}) \times U(\wt{\fg}) \rto & A\ot_k \Lambda(\fg)},
$$
defined by \hbox{$\wt{f}(u,v) = f(\pi(u),\pi(v))\ot_k 1$}, is a normal $2$-cocycle which satisfies the twisted module condition. Let
$$
\xymatrix{\eta\colon Y'_* \rto & (A\ot_k \Lambda(\fg))\#_{\wt{f}} U(\wt{\fg})}
$$
be the homomorphism of algebras defined by $\eta(a) = (a\ot_k 1)\# 1$ for all $a\in A$ and $\eta(y_x) = (1\ot_k 1)\# y_x$, $\eta(z_x) = (1\ot_k 1)\# z_x$ and $\eta(e_x) = (1\ot_k e_x)\# 1$ for all $x\in \fg$. Because of the Poincar\'e-Birkhoff-Witt theorem,
$$
\eta\bigl(y_{j_1}^{n_1}\cdots y_{j_h}^{n_h} \rho_{i_1}^{m_1}\cdots \rho_{i_l}^{m_l}\bigr)\qquad \left(h,l\ge 0 \text{,}\quad j_1<\dots<j_h \text{,}\atop i_1<\dots<i_l \in I\text{ and } m_j,n_j\ge 0 \right),
$$
is a basis of $(A\ot_k \Lambda(\fg))\#_{\wt{f}} U(\wt{\fg})$ as a left $A\ot_k \Lambda(\fg)$-module. The theorem follows easily from this fact.
\end{proof}

\begin{remark} A similar argument shows that each $Y_s$ is a free right $E$-module with the same basis.
\end{remark}

\begin{theorem}\label{th 1.2} Let $\wt{\mu}\colon Y_0\to E$ be the algebra map defined by $\wt{\mu}(a) = a$ for $a\in A$ and $\wt{\mu}(y_x) =\wt{\mu}(z_x) = 1\# x$ for $x\in \fg$. There is a unique derivation $\partial_*\colon Y_*\to Y_{*-1}$ such that $\partial(e_x) = z_x-y_x$ for $x\in\fg$. Moreover, the chain complex of $E$-bimodules
\[
\xymatrix{E&Y_0\lto_-{\wt{\mu}} &Y_1\lto_-{\partial_1} &Y_2\lto_-{\partial_2} &Y_3\lto_-{\partial_3} &Y_4\lto_-{\partial_4} &Y_5\lto_-{\partial_5}&\lto_-{\partial_6}\dots }
\]
is contractible as a complex of $(E,A)$-bimodules. A chain contracting homotopy
$$
\xymatrix{\si^{-1}_0\colon E\rto &Y_0},\qquad\xymatrix{\si^{-1}_{s+1}\colon Y_s\rto &Y_{s+1}}\quad (s\ge 0),
$$
is given by
\begin{align*}
&\si^{-1}(1) = 1,\\
&\si^{-1}\bigl(\rho_{i_1}^{m_1}e_{i_1}^{\de_1}\cdots \rho_{i_l}^{m_l} e_{i_l}^{\de_l} \bigr) = \begin{cases} \rho_{i_1}^{m_1-1}e_{i_1}\rho_{i_2}^{m_2}e_{i_2}^{\de_2}\cdots \rho_{i_l}^{m_l} e_{i_l}^{\de_l}& \text{if $\de_1 = 0$,}\\ 0 & \text{if $\de_1 = 1$,}\end{cases}
\end{align*}
where we assume that $i_1<\cdots<i_l$, $\de_1+\cdots+\de_l=s$ and $m_l+\de_l>0$.
\end{theorem}

\begin{proof} A direct computation shows that

\begin{itemize}

\smallskip

\item[-] $\wt{\mu}\xcirc \sigma^{-1}(1) = \wt{\mu}(1) = 1$,

\smallskip

\item[-] $\sigma^{-1}\xcirc\wt{\mu}(1) = \sigma^{-1}(1) = 1$ and $\partial\xcirc \sigma^{-1}(1) = \partial(0) = 0$,

\smallskip

\item[-] If $\bx = \rho_{i_1}^{m_1}\bx'$, where $m_1>0$ and $\bx' = \rho_{i_2}^{m_2}\cdots \rho_{i_l}^{m_l}$ with $i_1<\dots<i_l$, then
$$
\quad\sigma^{-1}\xcirc \wt{\mu}(\bx) = \sigma^{-1}(0) = 0\quad\text{and}\quad \partial \xcirc \sigma^{-1}(\bx) = \partial(\rho_{i_1}^{m_1-1}e_{i_1}\bx') = \bx,
$$

\smallskip

\item[-] Let $\bx = \rho_{i_1}^{m_1} e_{i_1}^{\delta_1}\bx'$, where $m_l+\delta_l>0$ and $\bx' = \rho_{i_2}^{m_2}e_{i_2}^{\delta_2}\cdots \rho_{i_l}^{m_l} e_{i_l}^{\delta_l}$ with $i_1<\dots<i_l$ and $\delta_1+\cdots+\delta_l = s > 0$. If $\delta_1 = 0$, then
\begin{align*}
\quad\qquad &\sigma^{-1}\xcirc\partial(\bx) = \sigma^{-1}\bigl(\rho_{i_1}^{m_1} \partial(\bx') \bigr) = \rho_{i_1}^{m_1-1} e_{i_1}\partial(\bx'),\\
& \partial\xcirc \sigma^{-1}(\bx) = \partial\bigl(\rho_{i_1}^{m_1-1} e_{i_1}\bx' \bigr) = \bx - \rho_{i_1}^{m_1-1} e_{i_1}\partial(\bx'),\\
\intertext{and if $\delta_1 = 1$, then}
& \sigma^{-1} \xcirc \partial(\bx) = \sigma^{-1}\bigl( \rho_{i_1}^{m_1+1} \bx' - \rho_{i_1}^{m_1}e_{i_1} \partial(\bx') \bigr) = \bx,\\
& \partial\xcirc \sigma^{-1}(\bx) = \partial(0) = 0.
\end{align*}
\end{itemize}
The result follows immediately from all these facts.
\end{proof}

\smallskip

For each $s\ge 0$ we consider $E\ot_k\fg^{\land s}$ as a right $K$-module via $(\bc\ot_k\bx) \lambda = \bc\lambda\ot_k\bx$. For $r,s\ge 0$, let $X_{rs} = (E\ot_k\fg^{\land s})\ot\ov{A}^r \ot E$. The groups $X_{rs}$ are $E$-bimodules in an obvious way. Let us consider the diagram of $E$-bimodules and $E$-bimodule maps
$$
\xymatrix{
\vdots\dto^{\partial_3}\\
Y_2\dto^{\partial_2} & X_{02}\lto_{\mu_2} & X_{12}\lto_{d^0_{12}}&\dots\lto_{d^0_{22}}\\
Y_1\dto^{\partial_1} & X_{01}\lto_{\mu_1} & X_{11}\lto_{d^0_{11}}&\dots\lto_{d^0_{21}}\\
Y_0 & X_{00}\lto_{\mu_0} & X_{10}\lto_{d^0_{10}}&\dots\lto_{d^0_{20}}, }
$$
where $\mu_*\colon X_{0*}\to Y_*$ and $d^0_{**}\colon X_{**}\to X_{*-1,*}$, are defined by:
\begin{align*}
&\mu(1\ot_k\bx_{1s}\ot 1) = e_{x_1}\dots e_{x_s},\\
\vspace{1.2\jot}
&d^0(1\ot_k\bx_{1s}\ot\ba_{1r}\ot 1) = (-1)^s a_1\ot_k\bx_{1s}\ot\ba_{2r}\ot 1\\
&\phantom{d^0(1\ot_k\bx_{1s}\ot\ba_{1r}\ot 1)} +\sum_{i=1}^{r-1} (-1)^{i+s}\ot_k\bx_{1s}\ot \ba_{1,i-1}\ot a_ia_{i+1}\ot\ba_{i+1,r}\ot 1\\
&\phantom{d^0(1\ot_k\bx_{1s}\ot\ba_{1r}\ot 1)} + (-1)^{r+s}\ot_k\bx_{1s}\ot\ba_{1,r-1}\ot a_r,
\end{align*}
Each horizontal complex in this diagram is contractible as a complex of $(E,K)$-bimodules. A chain contracting homotopy is the family
$$
\xymatrix{\si^0_{0s}\colon Y_s\rto & X_{0s}},\qquad \xymatrix{\si^0_{r+1,s}\colon X_{rs}\rto & X_{r+1,s}} \quad (r\ge 0),
$$
of $(E,K)$-bimodule maps, defined by
$$
\si^0(e_{x_1}\cdots e_{x_s} z_{x_{s+1}}\cdots z_{x_n}) =\sum_j a_j\ot_k\bx_{1s}\ot 1\# w_j,
$$
where $\sum_j a_j\# w_j = (1\# x_{s+1})\cdots (1\# x_n)$, and
$$
\si^0(1\ot_k\bx_{1s}\ot\ba_{1r}\ot a_{r+1}\# w) = (-1)^{r+s+1}\ot_k \bx_{1s}\ot \ba_{1,r+1}\ot 1\# w\quad(r\ge 0).
$$
(In order to see that the $\si^0$'s are right $K$-linear it is necessary to use that $K$ is stable under the action of $\fg$). Moreover, each $X_{rs}$ is a projective $E$-bimodule relative to the family of all epimorphism of $E$-bimodules which split as $(E,K)$-bimodule maps. We define $E$-bimodule maps
$$
\xymatrix{d^l_{rs}\colon X_{rs}\rto & X_{r+l-1,s-l}}\qquad\text{($r\ge 0$ and $1\le l\le s$),}
$$
recursively by:
$$
d^l(\mathbf{y}) = \begin{cases} -\si^0\xcirc\partial\xcirc\mu(\mathbf{y}) &\text{if $l=1$ and $r=0$,}\\
-\si^0\xcirc d^1\xcirc d^0 (\mathbf{y}) &\text{if $l=1$ and $r>0$,}\\
-\sum_{j=1}^{l-1}\si^0\xcirc d^{l-j}\xcirc d^j(\mathbf{y})&\text{if $l>1$ and $r=0$,}\\
-\sum_{j=0}^{l-1}\si^0\xcirc d^{l-j}\xcirc d^j(\mathbf{y}) &\text{if $l>1$ and $r>0$,}
\end{cases}
$$
where $\mathbf{y} = 1\ot_k\bx_{1s}\ot\ba_{1r}\ot 1$.

\begin{theorem}\label{th 1.3} The complex
\begin{equation}
\xymatrix{E&X_0\lto_-{\ov{\mu}} &X_1\lto_-{d_1} &X_2\lto_-{d_2} &X_3\lto_-{d_3} & X_4\lto_-{d_4} &X_5\lto_-{d_5}&\lto_-{d_6}\dots,}\label{f0}
\end{equation}
where
$$
\ov{\mu}(1\ot 1) = 1 ,\quad X_n =\bigoplus_{r+s=n} X_{rs}\quad\text{and}\quad d_n = \sum_{r+s=n\atop r+l> 0}\sum^s_{l=0} d^l_{rs},
$$
is a projective resolution of the $E$-bimodule $E$, relative to the family of all epimorphism of $E$-bimodules which split as $(E,K)$-bimodule maps. Moreover an explicit contracting
homotopy
$$
\xymatrix{\ov{\si}_0\colon E\rto & X_0},\qquad\xymatrix{\ov{\si}_{n+1}\colon X_n\rto & X_{n+1}}\quad (n\ge 0),
$$
of~\eqref{f0}, as a complex of $(E,K)$-bimodules, is given by
$$
\ov{\si}_0 = \si^0\xcirc\si^{-1}_0\quad\text{and}\quad\ov{\si}_{n+1} = -\sum_{l=0}^{n+1} \si^l_{l,n-l+1}\xcirc\si^{-1}_{n+1}\xcirc\mu_n + \sum_{r=0}^n\sum_{l=0}^{n-r} \si^l_{r+l+1,n-l-r},
$$
where
$$
\si^l_{l,s-l}\colon Y_s\to X_{l,s-l}\quad\text{and}\quad\si^l_{r+l+1,s-l}\colon X_{rs}\to X_{r+l+1,s-l}\quad(0<l\le s,r\ge 0)
$$
are recursively defined by
$$
\si^l = -\sum_{j=0}^{l-1}\si^0\xcirc d^{l-j}\xcirc\si^j.
$$
\end{theorem}

\begin{proof} By~\cite[Corollary A.2]{G-G2}.
\end{proof}

The boundary maps of the projective resolution of $E$ that we just found are defined recursively. Next we give closed formulas for them.

\begin{theorem}\label{th 1.4} For $x_i,x_j\in\fg$, we put $\wh{f}_{\! ij}= f(x_i,x_j) - f(x_j,x_i)$. We have:
\begin{align*}
& d^1(1\ot_k\bx_{1s}\ot\ba_{1r}\ot 1) =\sum_{i=1}^s (-1)^{i+1}\# x_i\ot_k \bx_{1\wh{\im}s} \ot \ba_{1r}\ot 1\\
&\phantom{d^1(1\ot_k\bx_{1s}\ot\ba_{1r}\ot 1)} +\sum_{i=1}^s (-1)^i\ot_k \bx_{1\wh{\im}s} \ot \ba_{1r}\ot 1\# x_i\\
&\phantom{d^1(1\ot_k\bx_{1s}\ot\ba_{1r}\ot 1)} +\sum_{i=1\atop 1\le h\le r}^s (-1)^i\ot_k \bx_{1\wh{\im}s}\ot\ba_{1,h-1}\ot a_h^{x_i}\ot\ba_{h+1,r}\ot 1\\
&\phantom{d^1(1\ot_k\bx_{1s}\ot\ba_{1r}\ot 1)} +\sum_{1\le i<j\le s} (-1)^{i+j}\ot_k [x_i,x_j]_{\fg}\land\bx_{1\wh{\im}\wh{\jm}s}\ot\ba_{1r}\ot 1,\\
\vspace{1.2\jot}
& d^2(1\ot_k\bx_{1s}\ot\ba_{1r}\ot 1) = \sum_{1\le i<j\le s\atop 0\le h\le r} (-1)^{i+j+h+s} \ot_k\bx_{1\wh{\im}\wh{\jm}s}\ot\ba_{1h}\ot\wh{f}_{\! ij}\ot\ba_{h+1,r}\ot 1
\end{align*}
and $d^l = 0$ for all $l\ge 3$.
\end{theorem}

\begin{proof} The proof of~\cite[Theorem 3.3]{G-G1} works in our more general context.
\end{proof}

\section{The comparison maps}
In this section we introduce and study comparison maps between $(X_*,d_*)$ and the canonical normalized Hochschild resolution $(E\ot\ov{E}^*\ot E,b'_*)$ of the $K$-algebra $E$. It is well known that there are morphisms of $E$-bimodule complexes
$$
\spreaddiagramcolumns{-0.5pc}
\xymatrix{\theta_*\colon (X_*,d_*)\rto & (E\ot\ov{E}^*\ot E,b'_*)} \quad \text{and} \quad
\xymatrix{\vartheta_*\colon (E\ot\ov{E}^*\ot E,b'_*)\rto & (X_*,d_*)},
$$
such that $\theta_0 =\vartheta_0 =\ide_{E\ot E}$ and that these morphisms are inverse one of each other up to homotopy. They can be recursively defined by $\theta_0 =\vartheta_0 = \ide_{E\ot E}$ and
\begin{align*}
&\theta(1\ot_k\bx_{1s}\ot\ba_{1r}\ot 1) = (-1)^n\theta\xcirc d(1\ot_k \bx_{1s}\ot\ba_{1r}\ot 1)\ot 1\\
\intertext{and}
&\vartheta(1\ot\bc_{1n}\ot 1) =\ov{\sigma}\xcirc\vartheta\xcirc b'(1\ot \bc_{1n}\ot 1),
\end{align*}
for $n\ge 1$, where $r+s = n$ and $\bc_{1n} = c_1\ot\cdots\ot c_n\in \ov{E}^n$. The following result was established without proof in~\cite{G-G1}.

\begin{proposition}\label{prop 2.1} We have:
$$
\theta\bigl(1\ot_k \bx_{1s}\ot \ba_{1r}\ot 1\bigr) = \sum_{\tau\in\mathfrak{S}_s} \sg(\tau)\ot \bigl(1\# x_{\tau(1)}\ot\cdots\ot 1\# x_{\tau(s)} \bigr)*\ba_{1r}\ot 1,
$$
where $\mathfrak{S}_s$ is the symmetric group in $s$ elements and $*$ denotes the shuffle product, which is defined by
$$
(\beta_1\ot\cdots\ot\beta_s)*(\beta_{s+1}\ot\cdots\ot\beta_n) =\sum_{\si\in\{(s,n-s)-\text{shuffles}\}} \sg(\si)\beta_{\si(1)}\ot\cdots\ot\beta_{\si(n)}.
$$
\end{proposition}

\begin{proof} We proceed by induction on $n=r+s$. The case $n=0$ is obvious. Suppose that $r+s = n$ and the result is valid for $\theta_{n-1}$. By the recursive definition of $\theta$, Theorem~\ref{th 1.4}, and the inductive hypothesis,
\begin{align*}
\theta(1\ot_k\bx_{1s}\ot\ba_{1r}\ot 1) & = (-1)^n\theta\xcirc d(1\ot_k\bx_{1s}\ot\ba_{1r}\ot 1) \ot 1\\
& = (-1)^n\theta\xcirc (d^0+d^1+d^2)(1\ot_k\bx_{1s}\ot\ba_{1r}\ot 1)\ot 1\\
&=\theta(1\ot_k\bx_{1s}\ot\ba_{1,r-1}\ot a_r)\ot 1\\
& +\theta\left(\sum_{i=1}^s (-1)^{i+n}\ot_k\bx_{1\wh{\im}s}\ot\ba_{1r}\ot 1\#x_i\right)\ot 1.
\end{align*}
The desired result follows now using again the inductive hypothesis.
\end{proof}

\begin{lemma}\label{le 2.2} Let $(g_i)_{i\in I}$ be the basis of $\fg$ considered in Theorem~\ref{th 1.1}. As in that theorem, let us write $e_i = e_{g_i}$ for each $i\in I$. The following facts hold:

\smallskip

\begin{enumerate}

\item $\ov{\sigma}_{n+1}\xcirc\ov{\sigma}_n = 0$ for all $n\ge 0$,

\smallskip

\item $\sigma^l((E\ot_k\fg^{\land s})\ot\ov{A}^r\ot K\# U(\fg)) = 0$ for all $0\le l\le s$,

\smallskip

\item $\sigma^l(e_{i_1}\cdots e_{i_n}) = 0$ for all $0<l\le n$,

\smallskip

\item $\sigma^l((E\ot_k\fg^{\land s})\ot\ov{A}^r\ot A) = 0$ for all $0<l\le s$,

\smallskip

\item $\sigma^{-1}\xcirc\mu\bigl(A\ot_k\fg^{\land n}\ot A\bigr) =0$,

\smallskip

\item Assume that $i_1<\cdots<i_n$. Then,
$$
\qquad\qquad\sigma^{-1}\xcirc\mu\bigl(1\ot_k g_{i_1}\land\cdots\land g_{i_n}\ot 1\# g_{i_{n+1}} \bigr) = \begin{cases} (-1)^ne_{i_1}\cdots e_{i_{n+1}} &\text{if $i_n<i_{n+1}$,}\\
0 &\text{otherwise.}\end{cases}
$$

\end{enumerate}

\end{lemma}

\begin{proof} (1)\enspace An inductive argument shows that there exist maps (which are left $E$-linear and right $K$-linear)
$$
\gamma^l_{rs}\colon X_{r+1,s}\to X_{r+l,s-l},
$$
such that $\sigma^l_{r+l+1,s-l} =\sigma^0_{r+l+1,s-l}\xcirc\gamma^l_{rs}\xcirc\sigma^0_{rs}$. Because of $\sigma^0\xcirc \sigma^0 = 0$, this implies that $\sigma^{l'}\xcirc \sigma^l = 0$, for all $l,l'\ge 0$. Thus,
$$
\ov{\sigma}_{n+1}\xcirc\ov{\sigma}_n = \sum_{l=0}^{n+1}\sigma^l\xcirc \sigma^{-1}\xcirc\mu \xcirc\sigma^0\xcirc\sigma^{-1}\xcirc\mu = 0,
$$
where the last equality holds because $\mu\xcirc\sigma^0 =\ide$ and $\sigma^{-1}\xcirc \sigma^{-1} = 0$.

\smallskip

\noindent (2)\enspace Since $\sigma^l =\sigma^0\xcirc\gamma^l\xcirc\sigma^0$ for $l>0$, we can assume that $l=0$. In this case the assertion follows immediately from the definition of $\sigma^0$.

\smallskip

\noindent (3)\enspace By the definition of $\si^0$ and Theorem~\ref{th 1.4},
\begin{align*}
&\sigma^0\xcirc d^1\xcirc\sigma^0(e_{i_1}\cdots e_{i_n}) =\sigma^0\xcirc d^1(1\ot_k g_{i_1} \land \cdots\land g_{i_n}\ot 1) = 0\\
\intertext{and}
&\sigma^0\xcirc d^2\xcirc\sigma^0(e_{i_1}\cdots e_{i_n}) =\sigma^0\xcirc d^2(1\ot_k g_{i_1} \land\cdots\land g_{i_n}\ot 1) = 0.
\end{align*}
Item~(3) follows now easily by induction on $l$, since, by the recursive definition of $\si^l$ and Theorem~\ref{th 1.4},
$$
\si^1 = -\si^0\xcirc d^1\xcirc \si^0\quad\text{and}\quad \si^l = -\si^0\xcirc d^1\xcirc \si^{l-1} -\si^0 \xcirc d^2\xcirc \si^{l-2}\,\text{ for $l\ge 2$.}
$$

\smallskip

\noindent (4)\enspace It is similar to the proof of item~(3).

\smallskip

\noindent (5)\enspace Since $e_ia=ae_i$ for all $i\in I$ and $a\in A$,
$$
\qquad\sigma^{-1}\xcirc\mu\bigl(a\ot_k g_{i_1}\land\cdots\land g_{i_n}\ot a'\bigr) = \sigma^{-1} \bigl(ae_{i_1}\cdots e_{i_n}a'\bigr) = \sigma^{-1}\bigl(aa'e_{i_1}\cdots e_{i_n}\bigr)= 0,
$$
where the last equality follows from the definition of $\sigma^{-1}$.

\smallskip

\noindent (6)\enspace We have
\begin{align*}
\qquad\quad \sigma^{-1}\xcirc\mu\bigl(1\ot_k g_{i_1}\land\cdots\land g_{i_n}\ot 1\# g_{i_{n+1}} \bigr) & = \sigma^{-1}\bigl(e_{i_1}\cdots e_{i_n} z_{i_{n+1}}\bigr)\\
& = \sigma^{-1}\bigl(e_{i_1}\cdots e_{i_n} (y_{i_{n+1}}+\rho_{i_{n+1}})\bigr)\\
& = \sigma^{-1}\bigl(y_{i_{n+1}} e_{i_1}\cdots e_{i_n}\bigr)\\
& + \sigma^{-1}\bigl(e_{i_1}\cdots e_{i_n} \rho_{i_{n+1}}\bigr),
\end{align*}
where $z_{i_{n+1}}$, $y_{i_{n+1}}$ and $\rho_{i_{n+1}}$ are as in Theorem~\ref{th 1.1}. So, in order to finish the proof it suffices to note that $\sigma^{-1}\bigl(y_{i_{n+1}}e_{i_1}\cdots e_{i_n}\bigr) = 0$ and
$$
\sigma^{-1}\bigl(e_{i_1}\cdots e_{i_n} \rho_{i_{n+1}} \bigr) = \begin{cases} (-1)^ne_{i_1} \cdots e_{i_{n+1}} &\text{if $i_n<i_{n+1}$,}\\ 0 &\text{otherwise,}\end{cases}
$$
which follows immediately from the fact that
$$
e_{i_j}\rho_{i_{n+1}} = \rho_{i_{n+1}}e_{i_j} +e_{[x_{i_j},x_{i_{n+1}}]_{\fg}}\quad\text{for all $j$ such that $i_j > i_{n+1}$,}
$$
and from the definition of $\sigma^{-1}$.
\end{proof}

\begin{theorem}\label{th 2.3} Let $(g_i)_{i\in I}$ be the basis of $\fg$ considered in Theorem~\ref{th 1.1}. Assume that $\bc_{1n} = c_1\ot\cdots\ot c_n \in\ov{E}^n$ is a simple tensor with $c_j\in A\cup\{1\# g_i : i\in I\}$ for all $j\in \{1,\dots,n\}$. If there exist $0\le s\le n$ and $i_1<\cdots<i_s$ in $I$, such that $c_j = 1\# g_{i_j}$ for $1\le j\le s$ and $c_j\in A$ for $s<j\le n$, then
$$
\vartheta(1\ot\bc_{1n}\ot 1) = 1\ot_k g_{i_1}\land\cdots\land g_{i_s}\ot\bc_{s+1,n}\ot 1.
$$
Otherwise, $\vartheta(1\ot\bc_{1n}\ot 1) = 0$.
\end{theorem}

\begin{proof} For all $n\ge 0$ we define $P_n$ by $\bc_{1n}\in P_n$ if there is $i_1<\cdots <i_s$ in $I$ such that $c_j = 1\# g_{i_j}$ for $j\le s$ and $c_j\in A$ for $j>s$. We now proceed by induction on $n$. The case $n=0$ is immediate. Assume that the result is valid for $\vartheta_n$. By item~(1) of Lemma~\ref{le 2.2} and the recursive definition of $\vartheta_n$, we have
$$
\ov{\sigma}\xcirc\vartheta(\bc'_{0n}\ot 1) =\ov{\sigma}\xcirc \ov{\sigma}\xcirc \vartheta \xcirc b'(\bc'_{0n}\ot 1) = 0,
$$
and so
$$
\vartheta(1\ot\bc_{1,n+1}\ot 1) = (-1)^{n+1} \ov{\sigma}\xcirc\vartheta(1\ot\bc_{1,n+1}).
$$
Assume that $c_j\in A\cup \{1\# g_i : i\in I\}$ for all $j\in\{1,\dots,n+1\}$. In order to finish the proof it suffices to show that

\begin{itemize}

\smallskip

\item[-] If $c_{1,n+1}\notin P_{n+1}$, then $\ov{\sigma}\xcirc\vartheta(1\ot\bc_{1,n+1}) = 0$,

\smallskip

\item[-] If $\bc_{1,n+1} = 1\# g_{i_1}\ot\cdots\ot 1\# g_{i_s}\ot\ba_{s+1,n+1}\in P_{n+1}$,
then
$$
\ov{\sigma}\xcirc\vartheta(1\ot\bc_{1,n+1}) =  (-1)^{n+1}\ot_k g_{i_1}\land \cdots\land g_{i_s}\ot\ba_{s+1,n+1}\ot 1.
$$

\smallskip

\end{itemize}
If $\bc_{1n}\notin P_n$, then $\vartheta(1\ot\bc_{1,n+1}) = 0$ by the inductive hypothesis. It remains to consider the case $\bc_{1n}\in P_n$. We divide this into three subcases.

\smallskip

\noindent 1) If $\bc_{1n} = 1\# g_{i_1}\ot\cdots\ot 1\# g_{i_s}\ot \ba_{s+1,n}$ and $c_{n+1} = a_{n+1}\in A$, then
\begin{align*}
\ov{\sigma}\xcirc\vartheta(1\ot\bc_{1,n+1}) & = \ov{\sigma}\bigl(1\ot_k g_{i_1}\land\cdots\land g_{i_s}\ot\ba_{s+1,n+1}\bigr)\\
& = \sigma^0\bigl(1\ot_k g_{i_1}\land\cdots\land g_{i_s}\ot\ba_{s+1,n+1}\bigr)\\
& = (-1)^{n+1}\ot_k g_{i_1}\land\cdots\land g_{i_s}\ot\ba_{s+1,n+1}\ot 1,
\end{align*}
by the inductive hypothesis, items~(4) and~(5) of Lemma~\ref{le 2.2}, and the definitions of $\ov{\sigma}$ and $\sigma^0$.

\smallskip

\noindent 2) If $\bc_{1n} = 1\# g_{i_1}\ot\cdots\ot 1\# g_{i_s}\ot \ba_{s+1,n}$ with $s<n$ and
$c_{n+1} = 1\# g_{i_{n+1}}$, then
$$
\ov{\sigma}\xcirc\vartheta(1\ot\bc_{1,n+1}) = \ov{\sigma}\bigl(1\ot_k g_{i_1}\land\cdots\land g_{i_s}\ot\ba_{s+1,n}\ot 1\# g_{i_{n+1}}\bigr) = 0,
$$
by the inductive hypothesis, the definition of $\ov{\sigma}$ and item~(2) of Lemma~\ref{le 2.2}.

\smallskip

\noindent 3) If $\bc_{1n} = 1\# g_{i_1}\ot\cdots\ot 1\# g_{i_n}$ and $c_{n+1} = 1\# g_{i_{n+1}}$, then
\begin{align*}
\ov{\sigma}\xcirc\vartheta(1\ot\bc_{1,n+1}) & = \ov{\sigma}\bigl(1\ot_k g_{i_1}\land\cdots\land g_{i_n}\ot 1\# g_{i_{n+1}}\bigr)\\
& = - \sigma^0\xcirc \sigma^{-1}\xcirc\mu\bigl(1\ot_k g_{i_1}\land\cdots\land g_{i_n}\ot 1\# g_{i_{n+1}}\bigr)\\
& = \begin{cases} (-1)^{n+1}\ot_k g_{i_1}\land \cdots\land g_{i_{n+1}}\ot 1 &\text{if $\bc_{1,n+1}\in P_{n+1}$,}\\0 &\text{otherwise,}\end{cases}
\end{align*}
by the inductive hypothesis, items~(2), (3) and~(6) of Lemma~\ref{le 2.2}, and the definitions of $\ov{\sigma}$ and $\sigma^0$.
\end{proof}

\section{The Hochschild cohomology}
Let $E= A\#_f U(\fg)$ and $M$ an $E$-bimodule. In this section we obtain a cochain complex $(\ov{X}_K^*(M),\ov{d}^*)$, simpler than the canonical one, giving the Hochschild cohomology of the $K$-algebra $E$ with coefficients in $M$.  When $K = k$ our result reduce to the one obtained in~\cite[Section 5]{G-G1}. Then, we obtain an expression that gives the cup product of the Hochschild cohomology of $E$ in terms of $(\ov{X}_K^*(E),\ov{d}^*)$. As usual, given $c\in E$ and $m\in M$, we let $[m,c]$ denote the commutator $mc-cm$.

\subsection{The complex $(\ov{X}_K^*(M),\ov{d}^*)$} For $r,s\ge 0$, let
$$
\ov{X}_K^{rs}(M) =\Hom_{K^e}(\ov{A}^r \ot_k \fg^{\land s},M),
$$
where $\ov{A}^r \ot_k \fg^{\land s}$ is considered as a $K$-bimodule via the canonical actions on $\ov{A}^r$. We define the morphism
$$
\spreaddiagramcolumns{-0.5pc}
\xymatrix{\ov{d}^{rs}_l\colon \ov{X}_K^{r+l-1,s-l}(M)\rto&\ov{X}_K^{rs}(M)}\quad \text{(with $0\le l\le\min(2,s)$ and $r+l>0$),}
$$
by:
\allowdisplaybreaks
\begin{align*}
&\ov{d}_0(\varphi)(\ba_{1r}\ot_k\bx_{1s}) = a_1\varphi(\ba_{2r}\ot_k\bx_{1s})\\
&\phantom{\ov{d}_0(\varphi)(\ba_{1r}\ot_k\bx_{1s})} + \sum_{i=1}^{r-1} (-1)^i\varphi(\ba_{1,i-1}\ot a_ia_{i+1}\ot\ba_{i+2,r}\ot_k\bx_{1s})\\
&\phantom{\ov{d}_0(\varphi)(\ba_{1r}\ot_k\bx_{1s})} + (-1)^r\varphi(\ba_{1,r-1} \ot_k \bx_{1s})a_r,\\
\vspace{1.5\jot}
&\ov{d}_1(\varphi)(\ba_{1r}\ot_k\bx_{1s}) = \sum_{i=1}^s (-1)^{i+r}\bigl[\varphi(\ba_{1r} \ot_k \bx_{1\wh{\im}s}),1\# x_i\bigr]\\
&\phantom{\ov{d}_1(\varphi)(\ba_{1r}\ot_k\bx_{1s})} +\sum_{i=1\atop 1\le h\le r}^s (-1)^{i+r} \varphi(\ba_{1,h-1}\ot a_h^{x_i}\ot\ba_{h+1,r}\ot_k\bx_{1\wh{\im}s})\\
&\phantom{\ov{d}_1(\varphi)(\ba_{1r}\ot_k\bx_{1s})} +\sum_{1\le i<j\le s} (-1)^{i+j+r} \varphi(\ba_{1r}\ot_k [x_i,x_j]_{\fg}\land\bx_{1\wh{\im}\wh{\jm}s})
\intertext{and}
&\ov{d}_2(\varphi)(\ba_{1r}\ot_k\bx_{1s}) =\sum_{1\le i<j\le s\atop 0\le h\le r} (-1)^{i+j+h} \varphi(\ba_{1h}\ot\wh{f}_{\! ij}\ot\ba_{h+1,r}\ot_k\bx_{1\wh{\im}\wh{\jm}s}),
\end{align*}
where $\wh{f}_{\! ij} = f(x_i,x_j) - f(x_j,x_i)$. Recall that $X_{rs} = (E\ot_k\fg^{\land s})\ot \ov{A}^r\ot E$. Applying the functor $\Hom_{E^e}(-,M)$ to the complex $(X_*,d_*)$ of Theorem~\ref{th 1.3}, and using Theorem~\ref{th 1.4} and the identifications $\gamma^{rs}\colon \ov{X}_K^{rs}(M)\to\Hom_{E^e}\bigl(X_{rs},M\bigr)$, given by
$$
\gamma(\varphi)(1\ot_k \bx_{1s}\ot\ba_{1r}\ot 1) = (-1)^{rs}\varphi(\ba_{1r}\ot_k \bx_{1s}),
$$
we obtain the complex
\[
%
\xymatrix{\ov{X}_K^0(M)\rto^-{\ov{d}^1} &\ov{X}_K^1(M)\rto^-{\ov{d}^2} & \ov{X}_K^2(M) \rto^-{\ov{d}^3} &\ov{X}_K^3(M)\rto^-{\ov{d}^4} &\ov{X}_K^4(M)\rto^-{\ov{d}^5} &\cdots,}
\]
where
$$
\ov{X}_K^n(M) = \bigoplus_{r+s=n} \ov{X}_K^{rs}(M)\qquad\hbox{and}\qquad \ov{d}^n = \sum_{r+s=n\atop r+l> 0}\sum_{l=0}^{\min(s,2)}\ov{d}_l^{rs}.
$$
Note that if $f(\fg\ot_k\fg )\subseteq K$, then $(\ov{X}_K^*(M),\ov{d}^*)$ is the total complex of the double complex $\bigl(\ov{X}_K^{**}(M),\ov{d}^{**}_0,\ov{d}^{**}_1\bigr)$.

\begin{theorem}\label{th 3.1} The Hochschild cohomology $\HS_K^*(E,M)$, of the $K$-algebra $E$ with coefficients in $M$, is the cohomology of $(\ov{X}_K^*(M),\ov{d}^*)$.
\end{theorem}

\begin{proof} It is an immediate consequence of the above discussion.
\end{proof}

\subsection{The comparison maps} The maps $\theta_*$ and $\vartheta_*$, introduced in Section~2, induce quasi-isomorphisms
$$
\xymatrix{\ov{\theta}^*\colon \bigl(\Hom_{K^e}\bigl(\ov{E}^*,M\bigr),b^*\bigr)\rto &(\ov{X}_K^*(M),\ov{d}^*)}
$$
and
$$
\xymatrix{\ov{\vartheta}^*\colon (\ov{X}_K^*(M),\ov{d}^*) \rto &\bigl(\Hom_{K^e} \bigl(\ov{E}^*,M\bigr), b^*\bigr)}
$$
which are inverse one of each other up to homotopy.

\begin{proposition}\label{prop 3.2} We have
$$
\ov{\theta}(\psi)(\ba_{1r}\ot_k \bx_{1s}) = \sum_{\tau\in\mathfrak{S}_s} (-1)^{rs} \sg(\tau)\psi\Bigr(\bigl(1\# x_{\tau(1)}\ot\cdots\ot 1\# x_{\tau(s)} \bigr)*\ba_{1r}\Bigl)
$$
\end{proposition}

\begin{proof} This follows immediately from Proposition~\ref{prop 2.1}.
\end{proof}

In the sequel we consider that $\ov{X}_K^{rs}\subseteq \ov{X}_K^{r+s}$ in the canonical way.

\begin{theorem}\label{th 3.3} Let $(g_i)_{i\in I}$ be the basis of $\fg$ considered in Theorem~\ref{th 1.1} and let $\varphi\in \ov{X}_K^{rs}$. Assume that $\bc_{1,r+s} = c_1\ot\cdots\ot c_{r+s} \in\ov{E}^{r+s}$ is a simple tensor with $c_j\in A\cup\{1\# g_i : i\in I\}$ for all $j\in \{1,\dots,r+s\}$. If $c_j = 1\# g_{i_j}$ with $i_1<\cdots<i_s$ in $I$ for $1\le j\le s$ and $c_j\in A$ for $s<j\le r+s$, then
$$
\ov{\vartheta}(\varphi)(\bc_{1,r+s}) = (-1)^{rs}\varphi(\bc_{s+1,r+s}\ot_k g_{i_1}\land\cdots\land g_{i_s}).
$$
Otherwise, $\ov{\vartheta}(\varphi)(\bc_{1,r+s}) = 0$.
\end{theorem}

\begin{proof} This follows immediately from Theorem~\ref{th 2.3}.
\end{proof}

As usual, in the following subsection we will write $\HH_K^*(E)$ instead of $\HS_K^*(E,E)$.

\subsection{The cup product} Recall that the cup product of $\HH_K^*(E)$ is given in terms of
$\bigl(\Hom_{K^e}\bigl(\ov{E}^*,E\bigr),b^*\bigr)$, by
\[
(\psi\smile \psi')(\bc_{1,m+n})= \psi(\bc_{1m}) \psi'(\bc_{m+1,m+n}),
\]
where $\psi\in \Hom_{K^e}(\ov{E}^m,E)$ and $\psi'\in\Hom_{K^e}(\ov{E}^n,E)$. In this subsection we compute the cup product in terms of the small complex $(\ov{X}_K^*(E),\ov{d}^*)$. Given $\varphi\in \ov{X}_K^{rs}(E)$ and $\varphi'\in \ov{X}_K^{r's'}(E)$ we define $\varphi\bullet \varphi'\in \ov{X}_K^{r+r',s+s'}(E)$ by
$$
(\varphi\bullet\varphi')(\ba_{1r''}\ot_k\bx_{1s''}) = \sum_{1\le j_1<\!\cdots<j_s\le s''} \sg(j_{1s}) \varphi(\ba_{1r}\ot_k \bx_{j_{1s}})\varphi'(\ba_{r+1,r''}\ot_k\bx_{l_{1s'}}),
$$
where

\begin{itemize}

\smallskip

\item[-] $\sg(j_{1s}) = (-1)^{\!r'\!s +\! \sum\limits_{u=1}^s\! (j_u\!-\!u)}$,

\smallskip

\item[-] $r'' = r+r'$ and $s'' = s+s'$,

\smallskip

\item[-] $1\le l_1<\cdots<l_{s'}\le s''$ denote the set defined by
$$
\{j_1,\dots,j_s\}\cup \{l_1,\dots,l_{s'}\} = \{1,\dots,s''\},
$$

\smallskip

\item[-] $\bx_{j_{1s}} = x_{j_1}\land\cdots \land x_{j_s}$ and $\bx_{l_{1s'}} = x_{l_1}\land\cdots \land x_{l_{s'}}$.

\end{itemize}

\begin{theorem}\label{th 3.4} The cup product of $\HH_K^*(E)$ is induced by the operation $\bullet$ in the complex $(\ov{X}_K^*(E),\ov{d}^*)$.
\end{theorem}

\begin{proof} Let $\varphi\in \ov{X}_K^{rs}(E)$ and $\varphi'\in \ov{X}_K^{r's'}(E)$. Let $r''$ and $s''$ be natural numbers satisfying $r''+s'' = r+r'+s+s'$ and let $\ba_{1r''}\ot_k\bx_{1s''}\in X^K_{r''s''}$. Let $(g_i)_{i\in I}$ be the basis of $\fg$ considered in Theorem~\ref{th 1.1}. Clearly we can assume that there exist $i_1<\cdots<i_{s''}$ in $I$ such that $x_j = g_{i_j}$ for all $1\le j\le s''$. By Proposition~\ref{prop 3.2},
$$
\ov{\theta}\bigl(\ov{\vartheta}(\varphi)\smile \ov{\vartheta}(\varphi')\bigr) (\ba_{1r''}\ot_k\bx_{1s''}) =  \bigl(\ov{\vartheta}(\varphi)\smile \ov{\vartheta}(\varphi')\bigr) (T)
$$
where
$$
T = \sum_{\tau\in \mathfrak{S}_{s''}} (-1)^{r''s''} \sg(\tau) \bigl((1\# x_{\tau(1)})\ot\cdots\ot (1\# x_{\tau(s'')}) \bigr)*\ba_{1r''}.
$$
In order to finish the proof it suffices to note that by Theorem~\ref{th 3.3}, this is zero if $r'' \ne r+r'$ and this is $(\varphi\bullet\varphi')(\ba_{1r''}\ot_k\bx_{1s''})$ if $r'' = r+r'$.
\end{proof}

\section{The Hochschild homology}
Let $E= A\#_f U(\fg)$ and $M$ an $E$-bimodule. In this section we obtain a chain complex $(\ov{X}^K_*(M),\ov{d}_*)$, simpler than the canonical one, giving the Hochschild homology of the $K$-algebra $E$ with coefficients in $M$. When $K = k$ our result reduce to the one obtained in~\cite[Section 4]{G-G1}. Then, we obtain an expression that gives the cap product of $\HS^K_*(E,M)$ in terms of $(\ov{X}_K^*(E),\ov{d}^*)$ and $(\ov{X}^K_*(E,M),\ov{d}_*)$. As in the previous section $[m,c]$ denotes the commutator $mc-cm$ of $m\in M$ and $c\in E$.

\subsection{The complex $(\ov{X}^K_*(M),\ov{d}_*)$} For $r,s\ge 0$, let
$$
\ov{X}^K_{rs}(M) = \frac{M\ot\ov{A}^r}{[M\ot\ov{A}^r,K]} \ot\fg^{\land s},
$$
where $[M\ot\ov{A}^r,K]$ is the $k$-vector space generated by the commutators $[m\ot\ba_{1r}, \lambda]$, with $\lambda\in K$ and $m\ot\ba_{1r}\in M\ot\ov{A}^r$. We define the morphism
$$
\spreaddiagramcolumns{-0.5pc}
\xymatrix{\ov{d}_{rs}^l\colon \ov{X}^K_{rs}(M)\rto &\ov{X}^K_{r+l-1,s-l}(M)}\quad \text{(with $0\le l\le\min(2,s)$ and $r+l>0$),}
$$
by:
\allowdisplaybreaks
\begin{align*}
&\ov{d}^0(\ov{m\ot\ba_{1r}}\ot_k \bx_{1s}) = \ov{ma_1\ot \ba_{2r}}\ot_k \bx_{1s}\\
&\phantom{\ov{d}^0(\ov{m\ot\ba_{1r}}\ot_k \bx_{1s})} + \sum_{i=1}^{r-1} (-1)^i \ov{m\ot \ba_{1,i-1} \ot a_ia_{i+1}\ot \ba_{i+2,r}}\ot_k \bx_{1s}\\
&\phantom{\ov{d}^0(\ov{m\ot\ba_{1r}}\ot_k \bx_{1s})}+ (-1)^r \ov{a_rm\ot \ba_{1,r-1}}\ot_k \bx_{1s},\\
\vspace{1.2\jot}
& \ov{d}^1(\ov{m\ot\ba_{1r}}\ot_k\bx_{1s}) = \sum_{i=1}^s (-1)^{i+r}\ov{[(1\# x_i),m]\ot \ba_{1r}}\ot_k\bx_{1\wh{\im}s} \\
&\phantom{d^1(\ov{m\ot\ba_{1r}}\ot_k\bx_{1s})} + \sum_{i=1\atop 1\le h\le r}^s(-1)^{i+r} \ov{m\ot\ba_{1,h-1}\ot a_h^{x_i}\ot\ba_{h+1,r}}\ot_k \bx_{1\wh{\im}s}\\
&\phantom{d^1(\ov{m\ot\ba_{1r}}\ot_k\bx_{1s})} + \sum_{1\le i<j\le s}(-1)^{i+j+r} \ov{m\ot \ba_{1r}}\ot_k [x_i,x_j]_{\fg}\land \bx_{1\wh{\im}\wh{\jm}s}
\intertext{and}
&\ov{d}^2(\ov{m\ot\ba_{1r}}\ot_k\bx_{1s}) = \sum_{1\le i<j\le s\atop 0\le h\le r} (-1)^{i+j+h} \ov{m\ot\ba_{1h}\ot\wh{f}_{\! ij}\ot\ba_{h+1,r}}\ot_k \bx_{1\wh{\im} \wh{\jm}s},
\end{align*}
where $\wh{f}_{\! ij}= f(x_i,x_j) - f(x_j,x_i)$ and $\ov{m\ot\ba_{1r}}$ denotes the class of $m\ot \ba_{1r}$ in $M\ot\ov{A}^r/[M\ot\ov{A}^r,K]$, etcetera. Recall that $X_{rs} = (E\ot_k\fg^{\land s})\ot\ov{A}^r\ot E$ and let $E^e$ be enveloping algebra of $E$. By tensoring on the left $X_{rs}$ over $E^e$ with $M$, and using Theorem~\ref{th 1.4} and the identifications $\gamma_{rs}\colon \ov{X}^K_{rs}(M)\to M \ot_{E^e} X_{rs}$, given by
$$
\gamma(\ov{m\ot\ba_{1r}}\ot_k \bx_{1s}) = (-1)^{rs} m\ot_{E^e}(1\ot_k\bx_{1s}\ot\ba_{1r}\ot 1),
$$
we obtain the complex
$$
%
\xymatrix{\ov{X}^K_0(M) &\ov{X}^K_1(M)\lto_-{\ov{d}_1} & \ov{X}^K_2(M) \lto_-{\ov{d}_2} & \ov{X}^K_3(M) \lto_-{\ov{d}_3} &\ov{X}^K_4(M)\lto_-{\ov{d}_4} &\cdots, \lto_-{\ov{d}_5}}
$$
where
$$
\ov{X}^K_n(M) =\bigoplus_{r+s=n} \ov{X}^K_{rs}(M)\qquad\hbox{and}\qquad \ov{d}_n = \sum_{r+s=n \atop r+l>0} \sum_{l=0}^{\min(s,2)}\ov{d}^l_{rs}.
$$
Note that if $f(\fg\ot_k\fg )\subseteq K$, then $(\ov{X}^K_*(M),\ov{d}_*)$ is the total complex of the double complex $\bigl(\ov{X}^K_{**}(M),\ov{d}_{**}^0,\ov{d}_{**}^1\bigr)$.

\begin{theorem}\label{th 4.1} The Hochschild homology $\HS^K_*(E,M)$, of the $K$-algebra $E$ with coefficients in $M$, is the homology of $(\ov{X}^K_*(M),\ov{d}_*)$.
\end{theorem}

\begin{proof} It is an immediate consequence of the above discussion.
\end{proof}

\subsection{The comparison maps} The maps $\theta_*$ and $\vartheta_*$, introduced in Section~2, induce quasi-isomorphisms
$$
\xymatrix{\ov{\theta}_*\colon (\ov{X}^K_*(M),\ov{d}_*)\rto &\left(\frac{M\ot \ov{E}^*} {[M\ot \ov{E}^*,K]},b_* \right)}
$$
and
$$
\xymatrix{\ov{\vartheta}_*\colon \left(\frac{M\ot\ov{E}^*}{[M\ot\ov{E}^*,K]},b_* \right) \rto &
(\ov{X}^K_*(M),\ov{d}_*)}
$$
which are inverse one of each other up to homotopy.

\begin{proposition}\label{prop 4.2} We have
$$
\ov{\theta}(\ov{m\ot \ba_{1r}}\ot_k\bx_{1s}) = \sum_{\tau\in\mathfrak{S}_s} (-1)^{rs} \sg(\tau) \ov{m\ot\bigl(1\# x_{\tau(1)}\ot\cdots\ot 1\# x_{\tau(s)} \bigr)*\ba_{1r}}
$$
\end{proposition}

\begin{proof} This follows immediately from Proposition~\ref{prop 2.1}.
\end{proof}

\begin{theorem}\label{th 4.3} Let $(g_i)_{i\in I}$ be the basis of $\fg$ considered in Theorem~\ref{th 1.1}. Assume that $\bc_{1n} = c_1\ot\cdots\ot c_n \in\ov{E}^n$ is a simple tensor with $c_j\in A\cup\{1\# g_i : i\in I\}$ for all $j\in \{1,\dots,n\}$. If there exist $0\le s\le n$ and $i_1<\cdots<i_s$ in $I$, such that $c_j = 1\# g_{i_j}$ for $1\le j\le s$ and $c_j\in A$ for $s<j\le n$, then
$$
\ov{\vartheta}(\ov{m\ot\bc_{1n}}) = (-1)^{s(n-s)}\ov{m\ot\bc_{s+1,n}}\ot_k g_{i_1}\land \cdots \land g_{i_s}.
$$
Otherwise, $\vartheta(\ov{m\ot\bc_{1n}}) = 0$.
\end{theorem}

\begin{proof} This follows immediately from Theorem~\ref{th 2.3}.
\end{proof}

\subsection{The cap product} Recall that the cap product
$$
\HS^K_p(E,M)\times \HH_K^q(E)\to \HS^K_{p-q}(E,M)\qquad (q\le p),
$$
is defined in terms of $\left(\frac{M\ot \ov{E}^*} {[M\ot \ov{E}^*,K]},b_* \right)$ and $\bigl(\Hom_{K^e}\bigl(\ov{E}^*,E\bigr),b^*\bigr)$, by
\[
\ov{m\ot\bc_{1p}}\smallfrown \psi = \ov{m \psi(\bc_{1q})\ot\bc_{q+1,p}},
\]
where $\psi\in \Hom_{K^e}(\ov{E}^q,E)$. In this subsection we compute the cup product in terms of the small complexes $(\ov{X}^K_*(M),\ov{d}_*)$ and $(\ov{X}_K^*(E),\ov{d}^*)$. Given
$$
\ov{m\ot \ba_{1r}}\ot_k\bx_{1s}\in \ov{X}^K_{rs}(M)\quad\text{and}\quad \varphi'\in \ov{X}_K^{r's'}(E)\qquad\text{with $r\ge r'$ and $s\ge s'$,}
$$
we define $(\ov{m\ot \ba_{1r}}\ot_k\bx_{1s})\bullet \varphi'\in \ov{X}^K_{r-r',s-s'}(M)$ by
$$
(\ov{m\ot \ba_{1r}}\ot_k\bx_{1s})\bullet\varphi'\! = \!\sum_{1\le j_1<\cdots<j_{s'}\le s}\!\!\! \sg(j_{1s'}) \ov{m\varphi'(\ba_{1r'}\!\ot_k\bx_{j_{1s'}}\!)\ot \ba_{r'+1,r}}\!\ot_k \! \bx_{l_{1,s-s'}}.
$$
where

\begin{itemize}

\smallskip

\item[-] $\sg(j_{1s'}) = (-1)^{\!r\!s'+\!r'\!s' +\! \sum\limits_{u=1}^{s'}\! (j_u\!-\!u)}$,

\smallskip

\item[-] $1\le l_1<\cdots<l_{s-s'}\le s$ denote the set defined by
$$
\{j_1,\dots,j_{s'}\}\cup \{l_1,\dots,l_{s-s'}\} = \{1,\dots,s\},
$$

\smallskip

\item[-] $\bx_{j_{1s'}} = x_{j_1}\land\cdots \land x_{j_{s'}}$ and $\bx_{l_{1,s-s'}} = x_{l_1}\land\cdots \land x_{l_{s-s'}}$.

\end{itemize}

\begin{theorem}\label{th 4.4} The cap product
$$
\HS^K_p(E,M)\times \HH_K^q(E)\to \HS^K_{p-q}(E,M),
$$
is induced by $\bullet$, in terms of the complexes $(\ov{X}^K_*(M),\ov{d}_*)$ and $(\ov{X}_K^*(E),\ov{d}^*)$.
\end{theorem}

\begin{proof} Let $\ov{m\ot \ba_{1r}}\ot_k\bx_{1s}\in \ov{X}^K_{rs}(M)$ and $\varphi'\in \ov{X}_K^{r's'}(E)$. Let $(g_i)_{i\in I}$ be the basis of $\fg$ considered in Theorem~\ref{th 1.1}. Clearly we can assume that there exist $i_1<\cdots<i_s$ in $I$ such that $x_j = g_{i_j}$ for all $1\le j\le s$. By Proposition~\ref{prop 4.2},
$$
\ov{\vartheta}\bigl(\ov{\theta}(\ov{m\ot \ba_{1r}}\ot_k\bx_{1s})\smallfrown \ov{\vartheta} (\varphi')\bigr) = \ov{\vartheta}\left(T \smallfrown \ov{\vartheta}(\varphi')\right),
$$
where
$$
T = \sum_{\sigma\in \mathfrak{S}_s} (-1)^{rs}\sg(\sigma)\bigl((1\# x_{\sigma(1)})\ot\cdots\ot (1\# x_{\sigma(s)}) \bigr)*\ba_{1r}.
$$
Hence, by Theorem~\ref{th 3.3}, if $r'>r$ or $s'>s$, then
$$
\ov{\vartheta}\bigl(\ov{\theta}(\ov{m\ot \ba_{1r}}\ot_k\bx_{1s})\smallfrown \ov{\vartheta}(\varphi')\bigr) = 0,
$$
and, if $r'\le r$ and $s'\le s$, then
$$
\ov{\vartheta}\bigl(\ov{\theta}(\ov{m\ot \ba_{1r}}\ot_k\bx_{1s})\smallfrown \ov{\vartheta}(\varphi')\bigr)\! = \! \sum_{1\le j_1<\cdots<j_{s'}\le s}\!\!\! \ov{\vartheta}\left(m\varphi'(\ba_{1r'}\! \ot_k\! \bx_{j_{1s'}})\ot T'_{l_{l,s-s'}}\!\right),
$$
where
$$
T'_{l_{l,s-s'}} = \sum_{\tau\in \mathfrak{S}_{s-s'}} (-1)^{rs+r's} \sg(\tau)\bigl((1\# x_{l_{\tau(1)}})\ot\cdots\ot (1\# x_{l_{\tau(s-s')}}) \bigr)*\ba_{r'+1,r}.
$$
In order to finish the proof it suffices to apply Theorem~\ref{th 4.3}.
\end{proof}

\section{The (co)homology of $S(V)\#_f U(\fg)$}
In this section we obtain a complexes $(\ov{Z}_*(M),\ov{\delta}_*)$ and $(\ov{Z}^*(M),\ov{\delta}^*)$, simpler than $(\ov{X}_K^*(M),\ov{d}^*)$ and $(\ov{X}^K_*(M),\ov{d}_*)$ respectively, giving the Hochschild homology of the $K$-algebra $E: = A\#_f U(\fg)$ with coefficients in an $E$-bimodule $M$

\begin{itemize}

\smallskip

\item[-] $K=k$ and $A$ is a symmetric algebra $S(V)$,

\smallskip

\item[-] $v^x\in k\oplus V$ for all $v\in V$ and $x\in \fg$,

\smallskip

\item[-] $f(x_1,x_2)\in k\oplus V$ for all $x_1,x_2\in \fg$.

\smallskip

\end{itemize}
Then, we obtain an expression that gives the cup product of $\HS_K^*(E,M)$ in terms of $(\ov{Z}^*(E),\ov{\delta}^*)$, and we obtain an expression that gives the cap product of $\HS^K_*(E,M)$ in terms of $(\ov{Z}_*(M),\ov{\delta}_*)$ and $(\ov{Z}^*(E),\ov{\delta}^*)$.

\smallskip

For $r,s\ge 0$, let $Z_{rs} = E\ot \fg^{\land s}\ot V^{\land r} \ot E$. The groups $Z_{rs}$ are $E$-bimodules in an obvious way. Let
$$
\xymatrix{\delta^l_{rs}\colon Z_{rs}\rto & Z_{r+l-1,s-l}}\qquad\text{($0\le l\le s$ and $r+l > 0$),}
$$
be the $E$-bimodule morphisms defined by
\begin{align*}
&\delta^0(1\ot\bx_{1s}\ot\bv_{1r}\ot 1) = \sum_{i=1}^r (-1)^{i+s} \bigl(v_i\ot \bx_{1s}\ot \bv_{1\wh{\im}r}\ot 1 - 1 \ot \bx_{1s}\ot \bv_{1\wh{\im}r}\ot v_i\bigr),\\
&\delta^1(1\ot\bx_{1s}\ot\bv_{1r}\ot 1) = \sum_{i=1}^s (-1)^{i+1}\# x_i\ot \bx_{1\wh{\im}s} \ot \bv_{1r}\ot 1\\
&\phantom{\delta^1(1\ot\bx_{1s}\ot\bv_{1r}\ot 1)} +\sum_{i=1}^s (-1)^i\ot \bx_{1\wh{\im}s} \ot \bv_{1r}\ot 1\# x_i\\
&\phantom{\delta^1(1\ot \bx_{1s}\ot\bv_{1r}\ot 1)} +\sum_{i=1\atop 1\le h\le r}^s (-1)^i\ot \bx_{1\wh{\im}s}\ot\bv_{1,h-1}\wedge v_h^{\ov{x}_i}\wedge \bv_{h+1,r}\ot 1\\
&\phantom{\delta^1(1\ot \bx_{1s}\ot\bv_{1r}\ot 1)} +\sum_{1\le i<j\le s} (-1)^{i+j}\ot [x_i,x_j]_{\fg}\land\bx_{1\wh{\im}\wh{\jm}s}\ot\bv_{1r}\ot 1
\intertext{and}
&\delta^2(1\ot\bx_{1s}\ot\bv_{1r}\ot 1) = \sum_{1\le i<j\le s} (-1)^{i+j+s} \ot\bx_{1\wh{\im} \wh{\jm}s}\ot\wh{f}_{\! ij}\wedge \bv_{1r}\ot 1,
\end{align*}
where $\bv_{hl} = v_h\wedge\cdots\wedge v_l$, $v_h^{\ov{x}_i}$ is the $V$-component of $v_h^{x_i}$ (that is $v_h^{\ov{x}_i}\in V$ and $v_h^{x_i}-v_h^{\ov{x}_i}\in k$) and $\wh{f}_{\! ij} = f_V(x_i,x_j) - f_V(x_j,x_i)$ in which $f_V(x_j,x_i)$ is the $V$-component of $f(x_j,x_i)$.

\begin{theorem}\label{th 5.1} The complex
\begin{equation*}
\xymatrix{E&Z_0\lto_-{\ov{\mu}} &Z_1\lto_-{\delta_1} &Z_2\lto_-{\delta_2} &Z_3\lto_-{\delta_3} &Z_4\lto_-{\delta_4} &Z_5\lto_-{\delta_5}&\lto_-{\delta_6}\dots,}
\end{equation*}
where
$$
\ov{\mu}(1\ot 1) = 1 ,\quad Z_n = \bigoplus_{r+s=n} Z_{rs}\quad\text{and}\quad \delta_n =
\sum_{r+s=n\atop r+l> 0}\sum^s_{l=0} \delta^l_{rs},
$$
is a projective resolution of the $E$-bimodule $E$. Moreover, the family of maps
$$
\Gamma_*\colon Z_*\to X_*,
$$
given by
$$
\Gamma(1\ot\bx_{1s}\ot\bv_{1r}\ot 1)=\sum_{\sigma\in \mathfrak{S}_r} \sg(\si)\ot\bx_{1s}\ot v_{\sigma(1)}\ot \cdots \ot v_{\sigma(r)}\ot 1,
$$
defines an morphism of $E$-bimodule complexes from $(Z_*,\delta_*)$ to $(X_*,d_*)$.
\end{theorem}

\begin{proof} It is clear that each $Z_n$ is a projective $E$-bimodule and a direct computation shows that $\Gamma_*$ is a morphism of complexes. Let
$$
G^0_* \subseteq G^1_* \subseteq G^2_* \subseteq G^3_*\subseteq \dots\quad\text{and}\quad F^0_* \subseteq F^1_* \subseteq F^2_* \subseteq F^3_*\subseteq \dots
$$
be the filtration of $Z_*$ and $X_*$ respectively, defined by
$$
G^i_n = \bigoplus_{r+s = n\atop s\le i} Z_{rs}\qquad\text{and}\qquad F^i_n =
\bigoplus_{r+s = n\atop s\le i} X_{rs}.
$$
In order to see that $\Gamma_*$ is a quasi-isomorphism it is sufficient to show that it induces a quasi-isomorphism between the graded complexes associated with the filtrations introduced above. In other words that the maps
$$
\xymatrix{\Gamma_{*s}\colon (Z_{*s},\delta^0_{*s}) \rto & (X_{*s},d^0_{*s})}\qquad (s\ge 0),
$$
defined by
$$
\Gamma(1\ot\bx_{1s}\ot\bv_{1r}\ot 1)=\sum_{\sigma\in \mathfrak{S}_r} \sg(\si)\ot\bx_{1s}\ot v_{\sigma(1)}\ot \cdots \ot v_{\sigma(r)}\ot 1,
$$
are quasi-isomorphisms, which follows easily from Proposition~\ref{prop 2.1}.
\end{proof}

\subsection{Hochschild cohomology} Let $M$ be an $E$-bimodule. For $r,s\ge 0$, let
$$
\ov{Z}^{rs}(M) = \Hom_k(V^r\ot \fg^{\land s},M).
$$
We define the morphism
$$
\spreaddiagramcolumns{-0.5pc}
\xymatrix{\ov{\delta}^{rs}_l\colon \ov{Z}^{r+l-1,s-l}(M)\rto&\ov{Z}^{rs}(M)}\quad \text{(with $0\le l\le\min(2,s)$ and $r+l>0$)}
$$
by:
\allowdisplaybreaks
\begin{align*}
&\ov{\delta}_0(\varphi)(\bv_{1r}\ot\bx_{1s}) = \sum_{i=1}^r (-1)^i [v_i, \varphi( \bv_{1\wh{\im}r}\ot  \bx_{1s})],\\
\vspace{1.5\jot}
&\ov{\delta}_1(\varphi)(\bv_{1r}\ot\bx_{1s}) = \sum_{i=1}^s (-1)^{i+r}\bigl[\varphi(\bv_{1r} \ot \bx_{1\wh{\im}s}),1\# x_i\bigr]\\
&\phantom{\ov{\delta}_1(\varphi)(\bv_{1r}\ot\bx_{1s})} +\sum_{i=1\atop 1\le h\le r}^s (-1)^{i+r} \varphi(\bv_{1,h-1}\wedge v_h^{\ov{x}_i}\wedge \bv_{h+1,r}\ot \bx_{1\wh{\im}s})\\
&\phantom{\ov{\delta}_1(\varphi)(\bv_{1r}\ot\bx_{1s})} +\sum_{1\le i<j\le s} (-1)^{i+j+r} \varphi(\bv_{1r}\ot [x_i,x_j]_{\fg}\land\bx_{1\wh{\im}\wh{\jm}s})
\intertext{and}
&\ov{\delta}_2(\varphi)(\bv_{1r}\ot\bx_{1s}) =\sum_{1\le i<j\le s} (-1)^{i+j} \varphi(\wh{f}_{\! ij} \wedge \bv_{1r}\ot\bx_{1\wh{\im} \wh{\jm}s}).
\end{align*}
Applying the functor $\Hom_{E^e}(-,M)$ to the complex $(Z_*,\delta_*)$, and using Theorem~\ref{th 5.1} and the identifications $\xi^{rs}\colon \ov{Z}^{rs}(M)\to \Hom_{E^e}\bigl(Z_{rs},M\bigr)$, given by
$$
\xi(\varphi)(1\ot \bx_{1s}\ot\bv_{1r}\ot 1) = (-1)^{rs}\varphi(\bv_{1r}\ot \bx_{1s}),
$$
we obtain the complex
\[
%
\xymatrix{\ov{Z}^0(M)\rto^-{\ov{\delta}^1} &\ov{Z}^1(M)\rto^-{\ov{\delta}^2} & \ov{Z}^2(M) \rto^-{\ov{\delta}^3} &\ov{Z}^3(M)\rto^-{\ov{\delta}^4} &\ov{Z}^4(M)\rto^-{\ov{\delta}^5} &\cdots,}
\]
where
$$
\ov{Z}^n(M) = \bigoplus_{r+s=n} \ov{Z}^{rs}(M)\qquad\hbox{and}\qquad \ov{\delta}^n = \sum_{r+s=n\atop r+l> 0}\sum_{l=0}^{\min(s,2)}\ov{\delta}_l^{rs}.
$$
Note that if $f(\fg\ot\fg )\subseteq k$, then $(\ov{Z}^*(M),\ov{\delta}^*)$ is the total
complex of the double complex $\bigl(\ov{Z}^{**}(M),\ov{\delta}^{**}_0,\ov{\delta}^{**}_1\bigr)$.

\begin{theorem}\label{th 5.2} The Hochschild cohomology $\HS^*(E,M)$, of $E$ with coefficients in $M$, is the cohomology of $(\ov{Z}^*(M),\ov{\delta}^*)$.
\end{theorem}

The map $\Gamma\colon (Z_*,\delta_*)\to (X_*,d_*)$ induces a quasi-isomorphism
$$
\xymatrix{\ov{\Gamma}^*\colon (\ov{X}_k^*(M),\ov{d}_*)\rto & (\ov{Z}^*(M),\ov{\delta}^*)}.
$$

\begin{proposition}\label{prop 5.3} We have
$$
\ov{\Gamma}(\varphi)(\bv_{1r}\ot\bx_{1s}) = \sum_{\sigma\in \mathfrak{S}_r} \sg(\si) \varphi( v_{\sigma(1)} \ot \cdots \ot v_{\sigma(r)}\ot \bx_{1s}).
$$
\end{proposition}

\begin{proof} This follows immediately from Theorem~\ref{th 5.1}.
\end{proof}

\subsection{The cup product} In this subsection we compute the cup product of $\HH^*(E)$ in terms of the complex $(\ov{Z}^*(E),\ov{\delta}^*)$. Given $\phi\in \ov{Z}^{rs}(E)$ and $\phi'\in \ov{Z}^{r's'}(E)$, we define $\phi\star \phi'\in \ov{Z}^{r+r',s+s'}(E)$ by
$$
(\phi\star \phi')(\bv_{1r''}\ot\bx_{1s''}) = \sum_{1\le i_1<\!\cdots<i_r\le r'' \atop 1\le j_1<\!\cdots<j_s\le s''}\sg(i_{\!1r},j_{\!1s}) \phi(\bv_{i_{1r}}\ot\bx_{j_{1s}}) \phi'(\bv_{h_{1r'}} \ot\bx_{l_{1s'}}),
$$
where

\begin{itemize}

\smallskip

\item[-] $\sg(i_{\!1r},j_{\!1s}) = (-1)^{\!r'\!s +\! \sum\limits_{u=1}^r\! (i_u\!-\!u) +\! \sum\limits_{u=1}^s\! (j_u\!-\!u)}$,

\smallskip

\item[-] $r'' = r+r'$ and $s'' = s+s'$,

\smallskip

\item[-] $1\le h_1<\cdots<h_{r'}\le r''$ denote the set defined by
$$
\{i_1,\dots,i_r\}\cup \{h_1,\dots,h_{r'}\} = \{1,\dots,r''\},
$$

\smallskip

\item[-] $1\le l_1<\cdots<l_{s'}\le s''$ denote the set defined by
$$
\{j_1,\dots,j_s\}\cup \{l_1,\dots,l_{s'}\} = \{1,\dots,s''\},
$$

\smallskip

\item[-] $\bv_{i_{1r}} = v_{i_1}\land\cdots \land v_{i_r}$ and $\bv_{h_{1r'}} = v_{h_1}\land\cdots \land v_{h_{r'}}$,

\smallskip

\item[-] $\bx_{j_{1s}} = x_{j_1}\land\cdots \land x_{j_s}$ and $\bx_{l_{1s'}} = x_{l_1}\land\cdots \land x_{l_{s'}}$.

\end{itemize}

\begin{theorem}\label{th 5.4} The cup product of $\HH^*(E)$ is induced by the operation $\star$ in the complex $(\ov{Z}^*(E),\ov{\delta}^*)$.
\end{theorem}

\begin{proof} By Theorem~\ref{th 3.4} it suffices to prove that
\begin{equation}
\ov{\Gamma}(\varphi\bullet \varphi') = \ov{\Gamma}(\varphi)\star \ov{\Gamma}(\varphi')\label{f1}
\end{equation}
for all $\varphi\in \ov{X}_k^{rs}(E)$ and $\varphi'\in \ov{X}_k^{r's'}(E)$. Let $\phi = \ov{\Gamma}(\varphi)$ and $\phi' = \ov{\Gamma}(\varphi')$. On one hand
\begin{align*}
(\phi\star \phi')& (\bv_{1r''}\ot \bx_{1s''}) = \sum_{1\le i_1<\!\cdots<i_r\le r'' \atop 1\le j_1<\!\cdots<j_s\le s''} \sg(i_{\!1r},j_{\!1s}) \phi(\bv_{i_{1r}}\ot\bx_{j_{1s}}) \phi'(\bv_{h_{1r'}} \ot\bx_{l_{1s'}})\\
& = \sum_{{1\le i_1<\!\cdots<i_r\le r'' \atop 1\le j_1<\!\cdots<j_s\le s''} \atop \tau\in \mathfrak{S}_r,\ \nu\in \mathfrak{S}_{r'}} \sg(i_{\!1r},j_{\!1s})\sg(\tau)\sg(\nu) \varphi (\bv_{i_{\tau(1r)}} \ot\bx_{j_{1s}}) \varphi'(\bv_{h_{\nu(1r')}} \ot\bx_{l_{1s'}}),
\end{align*}
where
$$
\bv_{i_{\tau(1r)}} = v_{i_{\tau(1)}}\ot \cdots\ot v_{i_{\tau(r)}}\quad\text{and}\quad \bv_{h_{\nu( 1r')}} = v_{h_{\nu(1)}}\ot \cdots\ot v_{h_{\nu(r')}}.
$$
On the other hand
\begin{align*}
\ov{\Gamma}(\varphi\bullet \varphi')&(\bv_{1r''}\ot \bx_{1s''}) = \sum_{\sigma \in \mathfrak{S}_{r''}}\sg(\sigma) (\varphi\bullet \varphi')(v_{\sigma(1)}\ot\cdots\ot v_{\sigma(r'')}\ot \bx_{1s})\\
& = \sum_{1\le j_1<\!\cdots<j_s\le s'' \atop \sigma \in \mathfrak{S}_{r''}} \sg(\sigma) \sg(j_{\!is})\varphi_1(v_{\sigma(1r)}\ot\bx_{j_{1s}})\varphi_2(v_{\sigma(r+1,r'')}\ot \bx_{l_{1s}}),
\end{align*}
where
$$
v_{\sigma(1r)} = v_{\sigma(1)}\ot\cdots\ot v_{\sigma(r)}\quad\text{and}\quad v_{\sigma(r+1,r'')} = v_{\sigma(r+1)}\ot\cdots \ot v_{\sigma(r'')}.
$$
Now, formula~\eqref{f1} if follows immediately from these facts.
\end{proof}

\subsection{Hochschild homology} Let $M$ be an $E$-bimodule. For $r,s\ge 0$, let
$$
\ov{Z}_{rs}(M) = M\ot V^{\land r} \ot \fg^{\land s}.
$$
We define the morphisms
$$
\xymatrix{\ov{\delta}^l_{rs}\colon \ov{Z}_{rs}(M)\rto & \ov{Z}_{r+l-1,s-l}(M)} \qquad\text{($0\le l\le s$ and $r+l > 0$)}
$$
by:
\begin{align*}
& \ov{\delta}^0(m\ot\bv_{1r}\ot\bx_{1s}) = \sum_{i=1}^r (-1)^i [m,v_i]\ot \bv_{1\wh{\im}r}\ot \bx_{1s},\\
\vspace{1.2\jot}
&\ov{\delta}^1(m\ot\bv_{1r}\ot\bx_{1s}) = \sum_{i=1}^s (-1)^{i+r} [1\# x_i,m]\ot \bv_{1r}\ot \bx_{1\wh{\im}s}\\
&\phantom{\delta^1(m\ot\bv_{1r}\ot \bx_{1s})} +\sum_{i=1\atop 1\le h\le r}^s (-1)^{i+r} m\ot\bv_{1,h-1}\wedge v_h^{\ov{x}_i}\wedge \bv_{h+1,r}\ot \bx_{1\wh{\im}s}\\
&\phantom{\delta^1(m\ot\bv_{1r}\ot \bx_{1s})} + \sum_{1\le i<j\le s} (-1)^{i+j+r} m\ot\bv_{1r}\ot [x_i,x_j]_{\fg}\land\bx_{1\wh{\im}\wh{\jm}s}
\intertext{and}
&\ov{\delta}^2(m\ot\bv_{1r}\ot\bx_{1s}) = \sum_{1\le i<j\le s} (-1)^{i+j} m\ot\wh{f}_{\! ij} \wedge \bv_{1r}\ot\bx_{1\wh{\im} \wh{\jm}s}.
\end{align*}
By tensoring on the left the complex $(Z_*,\delta_*)$ over $E^e$ with $M$, and using Theorem~\ref{th 5.1} and the identifications $\xi_{rs}\colon \ov{Z}_{rs}(M)\to M \ot_{E^e} Z_{rs}$, given by
$$
\xi(m\ot\bv_{1r}\ot \bx_{1s}) = (-1)^{rs} m\ot_{E^e}(1\ot\bx_{1s}\ot\bv_{1r}\ot 1),
$$
we obtain the complex
$$
%
\xymatrix{\ov{Z}_0(M) &\ov{Z}_1(M)\lto_-{\ov{\delta}_1} & \ov{Z}_2(M) \lto_-{\ov{\delta}_2} & \ov{Z}_3(M) \lto_-{\ov{\delta}_3} &\ov{Z}_4(M)\lto_-{\ov{\delta}_4} &\cdots, \lto_-{\ov{\delta}_5}}
$$
where
$$
\ov{Z}_n(M) =\bigoplus_{r+s=n} \ov{Z}_{rs}(M)\qquad\hbox{and}\qquad \ov{\delta}_n = \sum_{r+s=n\atop r+l>0} \sum_{l=0}^{\min(s,2)}\ov{\delta}^l_{rs}.
$$
Note that if $f(\fg\ot\fg )\subseteq k$, then $(\ov{Z}_*(M),\ov{\delta}_*)$ is the total complex of the double complex $\bigl(\ov{Z}_{**}(M),\ov{\delta}_{**}^0,\ov{\delta}_{**}^1 \bigr)$.

\begin{theorem}\label{th 5.5} The Hochschild homology $\HS_*(E,M)$, of $E$ with coefficients in $M$, is the homology of $(\ov{Z}_*(M),\ov{\delta}_*)$.
\end{theorem}

\begin{proof} It is an immediate consequence of the above discussion.
\end{proof}

The map $\Gamma\colon (Z_*,\delta_*)\to (X_*,d_*)$ induces a quasi-isomorphism
$$
\xymatrix{\ov{\Gamma}_*\colon (\ov{Z}_*(M),\ov{\delta}_*)\rto & (\ov{X}^k_*(M),\ov{d}_*)}.
$$

\begin{proposition}\label{prop 5.6} We have
$$
\ov{\Gamma}(m\ot \bv_{1r}\ot\bx_{1s}) = \sum_{\sigma\in \mathfrak{S}_r} \sg(\si)m\ot v_{\sigma(1)} \ot \cdots \ot v_{\sigma(r)}\ot \bx_{1s}.
$$
\end{proposition}

\begin{proof} This follows immediately from Theorem~\ref{th 5.1}.
\end{proof}

\subsection{The cap product} In this subsection we compute the cap product of
$$
\HS_p(E,M)\times \HH^q(E)\to \HS_{p-q}(E,M)\qquad (q\le p),
$$
in terms of the complexes $(\ov{Z}_*(M),\ov{\delta}_*)$ and $(\ov{Z}^*(E),\ov{\delta}^*)$. Given
$$
m\ot \bv_{1s}\ot\bx_{1s}\in \ov{Z}_{rs}(M)\quad\text{and}\quad \phi'\in \ov{Z}^{r's'}(E)\qquad\text{with $r\ge r'$ and $s\ge s'$,}
$$
we define $(m\ot \bv_{1r}\ot\bx_{1s})\star \phi'\in \ov{Z}_{r-r',s-s'}(M)$ by
$$
(m\ot \bv_{1r}\!\ot\bx_{1s})\star\phi'\! =\!\!\! \sum_{1\le i_1<\!\cdots<i_{r'}\le r \atop 1\le j_1<\!\cdots<j_{s'}\le s}\!\!\!\sg(i_{\!1r'},j_{\!1s'}) m\phi'(\bv_{i_{1r'}}\! \ot\bx_{j_{1s'}})\!\ot \bv_{h_{1,r'-r}}\!\ot\bx_{l_{1,s'-s}}.
$$
where

\begin{itemize}

\smallskip

\item[-] $\sg(i_{\!1r'},j_{\!1s'}) = (-1)^{\!r\!s' +\!r'\!s' +\! \sum\limits_{u=1}^{r'}\! (i_u\!-\!u) +\! \sum\limits_{u=1}^{s'}\! (j_u\!-\!u)}$,

\smallskip

\item[-] $1\le h_1<\cdots<h_{r-r'}\le r$ denote the set defined by
$$
\{i_1,\dots,i_{r'}\}\cup \{h_1,\dots,h_{r-r'}\} = \{1,\dots,r\},
$$

\smallskip

\item[-] $1\le l_1<\cdots<l_{s-s'}\le s$ denote the set defined by
$$
\{j_1,\dots,j_{s'}\}\cup \{l_1,\dots,l_{s-s'}\} = \{1,\dots,s\},
$$

\smallskip

\item[-] $\bv_{i_{1r'}} = v_{i_1}\land\cdots \land v_{i_{r'}}$ and $\bv_{h_{1,r-r'}} = v_{h_1}\land\cdots \land v_{h_{r-r'}}$,

\smallskip

\item[-] $\bx_{j_{1s'}} = x_{j_1}\land\cdots \land x_{j_{s'}}$ and $\bx_{l_{1,s-s'}} = x_{l_1}\land\cdots \land x_{l_{s-s'}}$.

\end{itemize}

\begin{theorem}\label{th 5.7} The cap product
$$
\HS_p(E,M)\times \HH^q(E)\to \HS_{p-q}(E,M)\qquad (q\le p),
$$
is induced by $\star$, in terms of the complexes $(\ov{Z}_*(M),\ov{\delta}_*)$ and $(\ov{Z}^*(E),\ov{\delta}^*)$.
\end{theorem}

\begin{proof} By Theorem~\ref{th 4.4} it suffices to prove that
\begin{equation}
\ov{\Gamma}(m\ot \bv_{1r}\ot\bx_{1s})\bullet \varphi' = \ov{\Gamma}\bigl((m\ot \bv_{1r}\ot\bx_{1s})\star \ov{\Gamma}(\varphi') \bigr)\label{f2}
\end{equation}
for all $m\ot \bv_{1r}\!\ot\bx_{1s}\in \ov{Z}_{rs}(M)$ and $\varphi'\in \ov{X}_k^{r's'}(E)$. Let $\phi' = \ov{\Gamma}(\varphi')$. On one hand
\begin{align*}
\ov{\Gamma}(m\ot \bv_{1r}\!\ot\bx_{1s}\!)\!\bullet\! \varphi' & =\!\!\!\! \sum_{1\le j_1<\!\cdots<j_{s'}\le s \atop \sigma \in \mathfrak{S}_r}\!\!\!\! \sg(\si j_{1s'}) m\varphi'(\bv_{\sigma(1r')}\!\ot\! \bx_{j_{1s'}})\!\ot v_{\sigma(r'+1,r)}\!\ot\! \bx_{l_{1,s-s'}},
\end{align*}
where $\sg(\si j_{1s'}) = \sg(\si)\sg(j_{1s'})$,
$$
v_{\sigma(1r')} = v_{\sigma(1)}\ot\cdots\ot v_{\sigma(r')}\quad\text{and}\quad v_{\sigma(r'+1,r)} = v_{\sigma(r'+1)}\ot\cdots \ot v_{\sigma(r)}.
$$
On the other hand
\begin{align*}
(m\ot \bv_{1r}&\ot\bx_{1s}) \star \phi' =\!\!\! \sum_{1\le i_1<\!\cdots<i_{r'}\le r \atop 1\le j_1<\!\cdots<j_{s'}\le s}\!\!\!\sg(i_{\!1r'},j_{\!1s'}) m\phi'(\bv_{i_{1r'}}\! \ot\bx_{j_{1s'}})\!\ot \bv_{h_{1,r'-r}}\!\ot\bx_{l_{1,s'-s}}\\
& = \! \sum_{{1\le i_1<\!\cdots<i_{r'}\le r \atop 1\le j_1<\!\cdots<j_{s'}\le s}\atop \tau \in \mathfrak{S}_{r'} } \!\!\!\sg(\tau)\sg(i_{\!1r'},j_{\!1s'}) m\varphi'(\bv_{i_{\tau(1r')}}\! \ot\bx_{j_{1s'}})\!\ot \bv_{h_{1,r'-r}}\!\ot\bx_{l_{1,s'-s}},
\end{align*}
where $\bv_{i_{\tau(1r')}} = v_{i_{\tau(1)}}\ot \cdots\ot v_{i_{\tau(r')}}$. Now, formula~\eqref{f2} if follows immediately from these facts.
\end{proof}


\begin{thebibliography}{B-G-R}

\bibitem[B-C-M]{B-C-M}
\newblock {R. J. Blattner, M. Cohen and S. Montgomery},
\newblock {\em Crossed products and inner actions of Hopf algebras},
\newblock {Trans. Amer. Math. Soc.}
\newblock {298}
\newblock {(1986)}
\newblock {671--711}.

\bibitem[B-G-R]{B-G-R}
\newblock {W. Borho, P. Gabriel and R. Rentschler},
\newblock {\em Primideale in Einh\"ullenden aufl\"osbarer Lie-Algebren},
\newblock {Lecture Notes in Mathematics}
\newblock {357},
\newblock {Springer-Verlag},
\newblock {Berlin Heidelberg New York}
\newblock {(1973)}.

\bibitem[Ch]{Ch}
\newblock {W. Chin},
\newblock {\em Prime ideals in differential operator rings and crossed products of infinite groups},
\newblock {Journal of Algebra}
\newblock {106}
\newblock {(1987)}
\newblock {78--104}.

\bibitem[D-T]{D-T}
\newblock {Y. Doi and M. Takeuchi},
\newblock {\em Cleft comodule algebras by a bialgebra},
\newblock {Comm. in Alg.}
\newblock {14}
\newblock {(1986)}
\newblock {801--817}.

\bibitem[G-G1]{G-G1}
\newblock {J.A. Guccione and J.J. Guccione},
\newblock {\em Hochschild (co)homology of differential opertators rings},
\newblock {Journal of Algebra}
\newblock {243}
\newblock {(2001)}
\newblock {596--614}.

\bibitem[G-G2]{G-G2}
\newblock {J.A. Guccione and J.J. Guccione},
\newblock {\em Hochschild (co)homology of a Hopf crossed products},
\newblock {K-theory}
\newblock {25}
\newblock {(2002)}
\newblock {139--169}.

\bibitem[Mc]{Mc}
\newblock {J. C. McConnell},
\newblock {\em Representations of solvable Lie algebras and the Gelfand-Kirillov conjecture},
\newblock {Proc. LMS}
\newblock {29}
\newblock {(1974)}
\newblock {453--484}.

\bibitem[Mc-R]{Mc-R}
\newblock {J. C. McConnell and J. C. Robson},
\newblock {\em No commutative Noetherian rings},
\newblock {Wiley-Interscience}
\newblock {New York}
\newblock {(1987)}.

\bibitem[M]{M}
\newblock {S. Montgomery},
\newblock {\em Crossed products of Hopf algebras and enveloping algebras},
\newblock{Perspectives in ring theory},
\newblock{Kluwer Academic},
\newblock{Dordrecht}
\newblock{(1988)}
\newblock{253--268}.

\bibitem[S]{S}
\newblock {R. Sridharan},
\newblock {\em Filtered algebras and representations of Lie algebras},
\newblock {Trans. Am. Math. Soc.}
\newblock {100}
\newblock {(1961)}
\newblock {530--550}.

\end{thebibliography}
\end{document}